\documentstyle[12pt]{article}
\begin{document}

\def\A{{\cal A}}
\def\B{{\cal B}}
\def\C{{\cal C}}
\def\D{{\cal D}}
\def\E{{\cal E}}
\def\F{{\cal F}}
\def\G{{\cal G}}
\def\H{{\cal H}}
\def\I{{\cal I}}
\def\J{{\cal J}}
\def\K{{\cal K}}
\def\L{{\cal L}}
\def\M{{\cal M}}
\def\N{{\cal N}}
\def\O{{\cal O}}
\def\P{{\cal P}}
\def\Q{{\cal Q}}
\def\R{{\cal R}}
\def\S{{\cal S}}
\def\T{{\cal T}}
\def\U{{\cal U}}
\def\V{{\cal V}}
\def\W{{\cal W}}
\def\X{{\cal X}}
\def\Y{{\cal Y}}
\def\Z{{\cal Z}}

\def\e{{\epsilon}}
\def\k{{\kappa}}
\def\l{{\ell}}

\def\argmax{\mathop{\rm arg\,max}}
\def\argmin{\mathop{\rm arg\,min}}
\def\qed{{\ \vrule width 1.5mm height 1.5mm \smallskip}}

\def\be{\begin{equation}}
\def\ee{\end{equation}}

% Equation Numbering 

\countdef\eqn=101
\eqn=0

\def\nr{\global\advance\eqn by 1
\eqno{(\the\eqn)}}

\def\0{\par\hangindent 10pt\noindent\hskip 10pt\hskip -10pt}

\countdef\refn=112
\refn=0

\def\Ref{\global\advance\refn by 1
\0 [\the \refn] }

\phantom{?}
\vskip 4cm

\centerline{\LARGE\bf 
Weakly Convergent Nonparametric Forecasting of}

\centerline{\LARGE\bf 
Stationary Time Series}

\vskip 4cm
\centerline{\LARGE\bf
Guszt\'av Morvai, 
Sidney Yakowitz and 
Paul Algoet}

\vskip 3 cm
\noindent
\centerline {\Large  IEEE Transactions on Information Theory Vol. 43, pp. 483-498, 1997.}
\vskip 3 cm

\begin{abstract}
The conditional distribution of the next outcome 
given the infinite past of a stationary process can be 
inferred from finite but growing segments of the past. 
Several schemes are known for constructing pointwise 
consistent estimates, but they all demand prohibitive 
amounts of input data. In this paper we consider 
real-valued time series and construct conditional 
distribution estimates that make much more efficient 
use of the input data. The estimates are consistent 
in a weak sense, and the question whether they are 
pointwise consistent is still open. For finite-alphabet 
processes one may rely on a universal data compression 
scheme like the Lempel-Ziv algorithm to construct 
conditional probability mass function estimates that 
are consistent in expected information divergence. 
Consistency in this strong sense cannot be attained 
in a universal sense for all stationary processes with 
values in an infinite alphabet, but weak consistency can. 
Some applications of the estimates to on-line forecasting, 
regression and classification are discussed.
\end{abstract}

\baselineskip 18 pt
\bigskip
\centerline
{\large\bf I. Introduction and Overview}

\medskip
We are motivated by some fundamental questions regarding 
inference of time series that were raised by T. Cover \cite{bib9} 
and concerning which significant progress has been made 
during the intervening years. The time series is a 
stationary process $\{X_t\}$ with values in a set $\X$
which may be a finite set, the real line, or a finite 
dimensional euclidean space. 
For $t\ge 0$ let $X^t=(X_0,X_1,\ldots,X_{t-1})$ denote the 
$t$-past at time $t$. It is also convenient to consider the 
outcome $X=X_0$, the $t$-past $X^{-t}=(X_{-t},\ldots,X_{-1})$ 
and the infinite past $X^-=(\ldots,X_{-2},X_{-1})$ at 
time $0$. The true process distribution $P$ is unknown 
a priori but is known to fall in the class $\P_s$ of 
stationary distributions on the sequence space $\X^\Z$.

Cover's list of questions included the following: 
given that $\{X_t\}$ is a $\{0,1\}$-valued time series 
with an unknown stationary ergodic distribution $P$, is it possible 
to infer estimates $\hat P\{X_t=1|X^t\}$ of the conditional 
probabilities $P\{X_t=1|X^t\}$ from the past $X^t$ such that 
\be[\hat P\{X_t=1|X^t\}-P\{X_t=1|X^t\}]\to 0 \quad\hbox{$P$-almost surely as $t\to\infty$?}\label{eq1}\ee
D. Bailey \cite{bib5} used the cutting and 
stacking technique of ergodic theory to prove that 
the answer is negative. A simple proof of this 
negative result is outlined in Proposition~3 of 
Ryabco \cite{bib30}. Bailey \cite{bib5} also discussed a result of 
Ornstein \cite{bib22} that provides a positive answer to 
a less demanding question of Cover \cite{bib9}, namely whether 
there exist estimates $\hat P\{X=1|X^{-t}\}$ based 
on the past $X^{-t}$ such that for all $P\in\P_s$,
\be\hat P\{X=1|X^{-t}\}\to P\{X=1|X^-\}
\quad\hbox{$P$-almost surely as $t\to\infty$.}\label{eq2}\ee
Ornstein constructed estimates $\hat P_k\{X=1|X^{-\lambda(k)}\}$
which depend on finite past segments 
$X^{-\lambda(k)}=(X_{-\lambda(k)},\ldots,X_{-1})$ and 
which converge almost surely to $P\{X=1|X^-\}$ for every 
$P\in\P_s$. The length $\lambda(k)$ of the data record 
$X^{-\lambda(k)}$ depends on the data itself, i.e. 
$\lambda(k)$ is a stopping time adapted to the filtration 
$\{\sigma(X^{-t}):\,t\ge 0\}$. To get estimates satisfying 
(\ref{eq2}), simply define $\hat P\{X=1|X^{-t}\}$ as the estimate 
$\hat P_k\{X=1|X^{-\lambda(k)}\}$ where $k$ is the largest 
integer such that $\hat P_k\{X=1|X^{-\lambda(k)}\}$ can be 
evaluated from the data $X^{-t}$ (that is, $X^{-\lambda(k)}$ 
is a suffix of the string $X^{-t}$ but $X^{-\lambda(k+1)}$ is not.) 
The true conditional probability $P\{X=1|X^{-t}\}$ converges 
to $P\{X=1|X^-\}$ almost surely by the martingale convergence 
theorem and the estimate $\hat P\{X=1|X^{-t}\}$ converges 
to the same limit, hence 
\be[\hat P\{X=1|X^{-t}\}-P\{X=1|X^{-t}\}]\to 0
\quad\hbox{$P$-almost surely and in $L^1(P)$.}\label{eq3}\ee
An on-line estimate $\hat P\{X_t=1|X^t\}$ can be 
constructed at time $t$ from the past $X^t$ in 
the same way as $\hat P\{X=1|X^{-t}\}$ was constructed 
from $X^{-t}$. By (\ref{eq3}) and stationarity 
\be[\hat P\{X_t=1|X^t\}-P\{X_t=1|X^t\}]\to 0
\quad\hbox{in $L^1(P)$ as $t\to\infty$.}\label{eq4}\ee
Thus the guessing scheme $\hat P\{X_t=1|X^t\}$ is
universally consistent in the weak sense of (\ref{eq4}), 
although no guessing scheme can be universally 
consistent in the pointwise sense of (\ref{eq1}).

Ornstein's result can be generalized when $\{X_t\}$ is 
a stationary process with values in a complete separable
metric (Polish) space $\X$. Algoet \cite{bib1} constructed estimates 
$\hat P_k(dx|X^{-\lambda(k)})$ that, with probability 
one under any $P\in\P_s$, converge in law to the true 
conditional distribution $P(dx|X^-)$ of $X=X_0$ given 
the infinite past. By setting 
$\hat P(dx|X^{-t})=\hat P_k(dx|X^{-\lambda(k)})$ for 
$\lambda(k)\le t<\lambda(k+1)$, one obtains estimates
$\hat P(dx|X^{-t})$ that almost surely converge in 
law to the random measure $P(dx|X^-)$ in the space of 
probability distributions on $\X$. Thus for any bounded 
continuous function $h(x)$ and any stationary
distribution $P\in\P_s$,
\be\int h(x)\,\hat P(dx|X^{-t})\to\int h(x)\,P(dx|X^-)
\quad\hbox{$P$-almost surely.}\label{eq5}\ee
A much simpler estimate $\hat P_k(dx|X^{-\lambda(k)})$ 
and convergence proof were obtained by Morvai, Yakowitz 
and Gy\"orfi \cite{bib21}. Their estimate $\hat P_k\{X\in B|X^{-\lambda(k)}\}$ 
of the conditional probability of a subset $B\subseteq\X$ 
has the structure of a sample mean:
\be\hat P_k\{X\in B|X^{-\lambda(k)}\}={1\over k}
\sum_{1\le i\le k}1\{X_{-\tau(i)}\in B\},\label{eq6}\ee
where the $X_{-\tau(i)}$ are samples of the process 
at selected instants in the past and $\lambda(k)$ 
is the smallest integer $t$ such that the indices 
$\{\tau(i):\,1\le i\le k\}$ can be inferred from the
segment $X^{-t}$. From careful reading of \cite{bib21}, one can 
surmise that $\lambda(k)$ will be huge for relatively 
small values of the sample size $k$. Morvai \cite{bib20} applied
the ergodic theorem for recurrence times of Ornstein and 
Weiss \cite{bib24} and argued that if $\{X_t\}$ is a stationary 
ergodic finite-alphabet process with positive entropy 
rate $H$ bits per symbol and $C$ is a constant such 
that $1<C<2^H$, then, with probability one, 
\be\lambda(k)\ge C^{C^{\cdot^{\cdot^C}}}
\quad\hbox{eventually for large $k$,}\label{eq7}\ee
where the height of the exponential tower is $k-k_0$ 
for some number $k_0$ that depends on the process
realization but not on $k$. To our knowledge, none 
of the strongly-consistent methods have been applied 
to any data sets, real or simulated.   

Scarpellini \cite{bib32} has applied the methods of Bailey \cite{bib5} 
and Ornstein \cite{bib22} to infer the conditional expectation 
$E\{X_\tau|\{X_s\}_{s\le 0}\}$ of the outcome $X_\tau$ 
at some fixed time $\tau>0$ given the infinite past 
of a stationary real-valued continuous-time process 
$\{X_t\}$ from past experience. The outcomes $X_t$ are 
assumed to be bounded in absolute value by some fixed 
constant $K$. Scarpellini constructs estimates by 
averaging samples taken at a finite number of regularly 
spaced instants in the past and proves that the 
estimates converge almost surely to the desired limit 
$E\{X_\tau|\{X_s\}_{s\le 0}\}$. His generalization 
of Ornstein's result is not quite straightforward, 
and the difficulty seems to be caused more by the 
continuity of the range space $[-K,K]$ than by the 
continuity of the time index $t$.

These works are of considerable theoretical interest 
because they point to the limits of what can be 
achieved by way of time series prediction. Pointwise 
consistency can be attained for all stationary processes, 
but the estimates are based on enormous data records. 
It is hard to say how much raw data are really needed 
to get estimates with reasonable precision. 
The nonparametric class of all stationary ergodic 
processes is very rich and can model all sorts of 
complex nonlinear dynamics with long range dependencies 
and periodicities at many different time scales. 
It is hopeless to get efficient estimates with bounds 
on the convergence rate unless one has a priori 
information that winnows the range of possibilities 
to some manageable subclass. 
In the literature on nonparametric estimation 
(e.g. see Gy\"orfi, H\"{a}rdle, Sarda and Vieu \cite{bib15} and also 
Marton and Shields \cite{bib19} 
), one imposes mixing 
conditions on the time series and then finds that the 
standard methods are consistent and achieve stated 
asymptotic rates of convergence. These approaches are 
preferable to the universal methods when one is assured 
of the mixing hypotheses. On the other hand, there is 
essentially no methodology for testing for mixing.

In the present study we relax the strong consistency 
requirement and push in the direction of greater
efficiency. Rather than demanding strong consistency 
or pointwise convergence in (\ref{eq5}), we shall be satisfied 
with weak consistency or mean convergence in $L^1(P)$.
(Note that mean convergence is equivalent to convergence 
in probability because the random variables are uniformly 
bounded.) Being more tolerant in this way enables us to 
significantly reduce the data demands of the algorithm. 
The estimates will again be defined as empirical averages 
of sample values, but the length of the raw data segment 
that must be inspected to collect a given number of  samples  
will grow only polynomially fast in the sample size (when $\X$ is a finite 
alphabet), rather than as a tower of exponentials in (\ref{eq7}).

For processes with values in a finite set $\X$, weak 
consistency means that for any stationary distribution 
$P$ on $\X^\Z$ and any $x\in\X$, the estimate 
$\hat P(x|X^{-t})=\hat P\{X=x|X^{-t}\}$ will converge in mean 
to the true conditional probability $P(x|X^-)=P\{X=x|X^-\}$:
\be\hat P(x|X^{-t})\to P(x|X^-)
\quad\hbox{in $L^1(P)$, for any $x\in\X$.}\label{eq8}\ee
There exist estimates that are universally consistent 
in a stronger sense. Given a universal data compression 
algorithm or a universal parsimonious modeling scheme 
for stationary processes with values in the finite
alphabet $\X$, we shall design estimates $\hat P(x|X^{-t})$
that are consistent in expected information divergence for 
all stationary $P$. The expectation of the Kullback-Leibler 
divergence between the conditional probability mass function 
$P(x|X^-)$ and the estimate $\hat P(x|X^{-t})$ will vanish 
in the limit as $t\to\infty$ for all $P\in\P_s$:
\be E_P\{I(P_{X|X^-}|\hat P_{X|X^{-t}})\}\to 0,\label{eq9}\ee
where
\be I(P_{X|X^-}|\hat P_{X|X^{-t}})=\sum_{x\in\X}
P(x|X^-)\log\left({P(x|X^-)\over \hat P(x|X^{-t})}\right).\label{eq10}\ee
Consistency in expected information divergence implies
consistency in mean as in (\ref{eq8}), and is equivalent to 
the requirement that for any stationary $P\in\P_s$
we have mean convergence
\be\log \hat P(X|X^{-t})\to \log P(X|X^-)\quad\hbox{in $L^1(P)$.}\label{eq11}\ee
The constructions of Ornstein \cite{bib22} and Morvai, Yakowitz and Gy\"{o}rfi  
\cite{bib21} yield estimates $\hat P(x|X^{-t})$ such that (\ref{eq11}) holds 
universally in the pointwise sense, but perhaps not in mean.

No estimates $\hat P(x|X^{-t})$ can be consistent in 
expected information divergence for all stationary 
processes with values in a countable infinite alphabet, 
but weak consistency as in (\ref{eq8}) is universally achievable. 
Barron, Gy\"orfi and van der Meulen \cite{bib7} consider an 
unknown distribution $P(dx)$ on an abstract measurable 
space $\X$ and construct estimates from independent 
samples so that the estimates are consistent in information 
divergence and in expected information divergence whenever 
$P(dx)$ has finite Kullback-Leibler divergence $I(P|M)<\infty$ 
relative to some known probability distribution $M(dx)$ on $\X$.  
In the present paper, the discussion of estimates that are consistent 
in expected information divergence is limited to 
the finite-alphabet case.

The organization of the paper is as follows.
In Section II we describe an algorithm for constructing
estimates $\hat P_k(dx|X^{-\lambda(k)})$ and prove weak 
consistency for all stationary real-valued time series. 
The method and its proof applies to time series with values 
in any $\sigma$-compact Polish space. In Section III 
we transform the estimates $\hat P_k(dx|X^{-\lambda(k)})$ 
into estimates $\hat P(dx|X^{-n})$ by letting $k$ depend 
on $n$. We choose an increasing sequence $k(n)$ and 
define the estimate $\hat P(dx|X^{-n})$ as 
$\hat P_{k(n)}(dx|X^{-\lambda(k(n))})$ if $\lambda(k(n))\le n$ 
and as some default measure $Q(dx)$ otherwise. 
If $k(n)$ grows sufficiently slowly with $n$ then 
the data requirement $\lambda(k(n))$ will seldom 
exceed the available length $n$ and the estimates 
$\hat P(dx|X^{-n})$ will be weakly consistent just like 
the estimates $\hat P_{k(n)}(dx|X^{-\lambda(k(n))})$.
Section IV is about modeling and data compression 
and about estimates that are consistent in expected 
information divergence for stationary processes 
with values in a finite alphabet. In Section V, 
we shift $\hat P(dx|X^{-t})$ from time $0$ to 
time $t$ and show that the shifted estimates 
$\hat P(dx_t|X^t)$ can be used for sequential 
forecasting or on-line prediction. We show that 
one can make sequential decisions based on the 
shifted estimates $\hat P(dx_t|X^t)$ so that the 
average loss per decision converges in mean to the 
minimum long run average loss that could be attained 
if one could make decisions with knowledge of the 
true conditional distribution of the next outcome
given the infinite past at each step. In particular, 
the average rate of incorrect guesses in classification 
and the average of the mean squared error in regression 
converge to the minimum that could be attained 
if the infinite past were known to begin with.

We would like to alert the reader about some of our 
notational conventions. Only one level of subscripts 
or superscripts is allowed in equations that are 
embedded in the text and so we are often forced to 
adopt the flat functional notation $\lambda(k)$, 
$\lambda(k(n))$, $\l(k)$, $J(k)$, $\tau(k,j)$, etc. 
However, the equations sometimes look better with nested 
subscripts and superscripts and therefore we prefer 
to write $\lambda_k$, $\lambda_{k(n)}$, $\l_k$, $J_k$, 
$\tau^k_j$, etc.  in the displayed equations. 
We hope that mixing of these notational conventions 
will not be a source of confusion but rather will 
improve the readability of the paper. Logarithms 
and entropy rates are taken in base $2$ unless 
specified otherwise, and exponential growth rates 
are really doubling rates. 

\bigskip
\centerline{\large\bf II. Learning the Conditional 
Distribution $P(dx|X^-)$}

\medskip
Let $\{X_t\}$ be a real-valued stationary time series.
The process distribution is unknown but shift-invariant.
We wish to infer the conditional distribution of $X=X_0$ 
given the infinite past $X^-$ from past experience. 
We show that it is very easy to construct weakly 
consistent estimates $\hat P_k(dx|X^{-\lambda(k)})$ 
depending on finite past data segments $X^{-\lambda(k)}$
such that for every bounded continuous function $h(x)$ 
on $\X$ and any stationary distribution $P\in\P_s$,
\be\lim_k\,\int h(x)\hat P_k(dx|X^{-\lambda(k)})=
\int h(x)\,P(dx|X^-)\quad\hbox{in $L^1(P)$.}\label{eq12}\ee

The estimates $\hat P_k(dx|X^{-\lambda(k)})$ will
be defined in terms of quantized versions of the process 
$\{X_t\}$. Let $\X$ denote the real line and let 
$\{\B_k\}_{k\ge 1}$ be an increasing sequence of 
finite subfields that asymptotically generate the 
Borel $\sigma$-field on $\X$. Let $x\mapsto [x]^k$ 
denote the quantizer that maps any point $x\in\X$ 
to the atom of $\B_k$ that happens to contain $x$. 
For any integer $\l\ge 1$ let $[X^{-\l}]^k$ denote the 
quantized sequence $([X_{-\l}]^k,\ldots,[X_{-1}]^k)$. 
Given any integer $J\ge 1$, one may search backwards
in time and collect $J$ samples of the process 
at times when the quantized $\l$-past looks 
exactly like the quantized $\l$-past at time $0$.
Let $\lambda=\lambda(k,\l,J)$ denote the length of 
the data segment $X^{-\lambda}=(X_{-\lambda},\ldots,X_{-1})$ 
that must be inspected to find these $J$ samples 
and let $\hat P_{k,\l,J}(dx|X^{-\lambda})$ denote 
the empirical distribution of those samples.
Then $\hat P_{k,\l,J}(dx|X^{-\lambda})$ will be 
a good estimate of $P(dx|X^-)$ if the sample size $J$, 
the context length $\l$ and the quantizer index $k$ 
are sufficiently large. In fact, if $k$ and $\l$ are
fixed and the sample size $J$ tends to infinity then
by the ergodic theorem, $\hat P_{k,\l,J}(dx|X^{-\lambda(k,\l,J)})$
will converge in law to $P(dx|[X^{-\l}]^k)$. If we now 
refine the context by increasing $k$ and $\l$, then 
$P(dx|[X^{-\l}]^k)$ will converge in law to $P(dx|X^-)$ 
by the martingale convergence theorem. The question is 
how to turn this limit of limits into a single limit by 
letting $k,\l$ and $J$ increase simultaneously to infinity.
We must make $k$ and $\l$ large to reduce the bias 
and we must make $J$ large to reduce the variance of 
the estimates. We will let $\l$ and $J$ grow with $k$ 
and show that if $\l(k)$ and $J(k)$ are monotonically 
increasing to infinity then the empirical conditional 
distribution estimate $\hat P_k(dx|X^{-\lambda(k)})=
\hat P_{k,\l(k),J(k)}(dx|X^{-\lambda(k,\l(k),J(k))})$
converges weakly to $P(dx|X^-)$. After this brief 
outline we now proceed with a detailed development.

Let $\{\l_k\}_{k\ge 1}$ and $\{J_k\}_{k\ge 1}$ be two 
nondecreasing unbounded sequences of positive integers.
We often write $\l(k)$ and $J(k)$ instead of 
$\l_k$ and $J_k$. For fixed $k\ge 1$ let 
$\{-\tau^k_j\}_{j\ge 0}$ and $\{\tilde\tau^k_j\}_{j\ge 0}$ 
denote the sequences of past and future recurrence 
times of the pattern $[X^{-\l(k)}]^k$. Thus we set 
$\tau^k_0=\tilde\tau^k_0=0$ and for $j=1,2,\ldots$ 
we inductively define 
\be\tau^k_j=\min\,\{t>\tau^k_{j-1}:\,
([X_{-\l_k-t}]^k,\ldots,[X_{-1-t}]^k)=
([X_{-\l_k}]^k,\ldots,[X_{-1}]^k)\},\label{eq13}\ee
\be\tilde \tau^k_j=\min\,\{t>\tilde\tau^k_{j-1}:\,
([X_{-\l_k}]^k,\ldots,[X_{-1}]^k)
=([X_{-\l_k+t}]^k,\ldots,[X_{-1+t}]^k)\}.\label{eq14}\ee
The random variables $\tau(k,j)=\tau^k_j$ and 
$\tilde\tau(k,j)=\tilde\tau^k_j$ are finite almost 
surely by Poincar\'e's recurrence theorem for the 
quantized process $\{[X_t]^k\}$, cf. Theorem~6.4.1 
of Gray \cite{bib14}. The lengths $\lambda_k=\lambda(k)$ 
and estimates $\hat P_k(dx|X^{-\lambda(k)})$ are 
now defined by the formulas
\be\lambda_k=\lambda(k)=\l(k)+\tau(k,J_k),\label{eq15}\ee
\be\hat P_k(dx|X^{-\lambda_k})=
{1\over J_k}\sum_{1\le j\le J_k}\delta_{X_{-\tau(k,j)}}(dx),\label{eq16}\ee
where $\delta_\xi(dx)$ is the Dirac measure that places 
unit mass at the point $\xi\in\X$. Thus for any Borel
set $B$, the conditional probability estimate
\be\hat P_k\{X\in B|X^{-\lambda_k}\}=
{1\over J_k}\sum_{1\le j\le J_k}1\{X_{-\tau(k,j)}\in B\}\label{eq17}\ee
is obtained by searching for the $J_k$ most recent 
occurrences of the pattern $[X^{-\l(k)}]^k$ and calculating 
the relative frequency with which the next realized 
symbols $X_{-\tau(k,j)}$ hit the set $B$. We shall 
prove that $\hat P_k(dx|X^{-\lambda(k)})$ is a weakly 
consistent estimate of $P(dx|X^-)$. The precise 
statement and the proof are broken down in two parts.

\smallskip
\noindent
{\bf Theorem 1A.} {\sl For any set $B$ in the generating 
field $\bigcup_k\,\B_k$ and any stationary process 
distribution $P\in\P_s$ we have mean convergence
\be\lim_k\hat P_k\{X\in B|X^{-\lambda_k}\}
=P\{X\in B|X^-\}\quad\hbox{in $L^1(P)$.}\label{eq18}\ee
}

The proof is somewhat technical and is placed in 
the Appendix. In the second part we argue that 
the estimators $\hat P_k(dx|X^{-\lambda(k)})$ 
can be employed to infer the regression function 
$E\{h(X)|X^-\}=\int h(x)\,P(dx|X^-)$ of any 
bounded continuous function $h(x)$ given the past.

\smallskip
\noindent
{\bf Theorem 1B.}
{\sl Let $\{X_t\}$ be a real-valued stationary time series. 
If the fields $\B_k$ are generated by intervals and the 
estimator $\hat P_k(dx|X^{-\lambda(k)})$ is defined as in (\ref{eq16}) 
then for any bounded continuous function $h(x)$ on $\X$,
\be\lim_k\,\int h(x)\hat P_k(dx|X^{-\lambda_k})=
\int h(x)\,P(dx|X^-)\quad\hbox{in $L^1(P)$.}\label{eq19}\ee 
}

\noindent
{\sl Proof:} 
Pick some bound $M$ such that $|h(x)|\le M$ on $\X$. 
Given $\e>0$ there exists an integer $\k$ and 
a finite interval $K$ in the field $\B_\k$ such that 
\be P\{X\in K\}>1-{\e\over M}.\label{eq20}\ee
If necessary we increase $\k$ until $\k$ is sufficiently
large so that there exists a $\B_\k$-measurable function 
$g(x)$ such that $|h(x)-g(x)|\le \e$ on $K$. Assuming 
$g(x)=0$ outside $K$, we have
\be|h(x)-g(x)|\le f(x)=\e\, 1\{x\in K\}+M \,1\{x\not\in K\}.\label{eq21}\ee
Let $\hat P_k$ and $P^-$ be shorthand for 
$\hat P_k(dx|X^{-\lambda(k)})$ and $P(dx|X^-)$. Then
\be\bigg|\int h\,d\hat P_k-\int h\,dP^-\bigg|
\le \int|h-g|\,d\hat P_k+
\bigg|\int g\,d\hat P_k-\int g\,dP^-\bigg|
+\int|g-h|\,dP^-.\label{eq22}\ee
The function $g(x)$ is a  finite linear combination of 
indicator functions of $\B_\k$-measurable subsets,
and Theorem 1A implies that $\int g\,d\hat P_k$ 
converges to $\int g\,dP^-$ in $L^1$:
\be E\bigg|\int g\,d\hat P_k-\int g\,dP^-\bigg|\to 0.\label{eq23}\ee 
The function $f(x)$ is $\B_\k$-measurable and bounded,
hence $\int f\,d\hat P_k$ converges to $\int f\,dP^-$ 
in $L^1$ and the expectations converge:
\be E\int f\,d\hat P_k\to E\int f\,dP^-=Ef.\label{eq24}\ee
Since $|h-g|\le f$ and
$Ef\le \e\,P\{X\in K\}+M\,P\{X\not\in K\}<2\e$ by (\ref{eq20}) 
and (\ref{eq21}), it follows from (\ref{eq22}), (\ref{eq23}) and (\ref{eq24}) that
\be E\bigg|\int h\,d\hat P_k-\int h\,dP^-\bigg|
\le 2\e+\e+2\e\quad\hbox{eventually for large $k$.}\label{eq25}\ee
Thus $E|\int h\,d\hat P_k-\int h\,dP^-|\to 0$,
and this is the desired conclusion (\ref{eq19}).
\qed

\smallskip
Theorem 1B holds in general if $\X$ is a $\sigma$-compact 
Polish space and the fields $\B_k$ are suitably chosen. 
Indeed, let $\{K_k\}_{k\ge 1}$ be an increasing sequence 
of compact subsets with union $\bigcup_k\,K_k=\X$. 
For any fixed $k$ one may cover $K_k$ with a finite 
collection of open balls having diameter less than 
$\e_k$, where $\e_k\searrow 0$ as $k\to\infty$.    
Let $\B_k$ denote the smallest field containing
$\B_{k-1}$ and the sets $B\cap K_k$ where $B$ ranges
over all balls in the finite cover of $K_k$. (We 
start with the trivial field $\B_0=\{\emptyset,\X\}$.)
Any bounded continuous function $h(x)$ on $\X$ is
uniformly continuous on each compact subset of $\X$.
If $|h(x)|\le M$ and $\e>0$, then for sufficiently
large $\k$ there exists some compact subset $K$ in 
$\B_\k$ such that $P\{X\not\in K\}\le \e/M$ and 
$h(x)$ oscillates less than $\e$ on each atom of 
$\B_\k$ that is contained in $K$. Thus there exists 
a $\B_\k$-measurable function $g(x)$ such that 
$|h(x)-g(x)|<\e$ on $K$ and $g(x)=0$ outside $K$.
We can then proceed as in the proof of Theorem~1B   
to prove that for any bounded continuous 
function $h(x)$,
\be\int h(x)\,\hat P_k(dx|X^{-\lambda(k)})\to E\{h(X)|X^-\} 
\quad\hbox{in $L^1$.}\label{eq26}\ee

\bigskip
\centerline
{\large\bf III. Truncation of the Search Depth }

\medskip
The estimates $\hat P_k(dx|X^{-\lambda(k)})$ are 
based on finite but random length segments of the past. 
We shall transform these into estimates $\hat P(dx|X^{-n})$ 
that depend on finite past segments with deterministic 
length but that still are weakly consistent. The details 
are somewhat more involved than for the strongly consistent 
estimates in Section I. In terms of the empirical conditional 
distribution $\hat P_{k,\l,J}(dx|X^{-\lambda(k,\l,J)})$ that 
was defined in the outline of Section II, the question is 
how fast $k$, $\l$ and $J$ may increase with $n$ so that 
$\lambda(k(n),\l(n),J(n))\le n$ with high probability.
The weak consistency of the estimates 
$\hat P_{k(n),\l(n),J(n)}(dx|X^{-\lambda(k(n),\l(n),J(n))})$
will not suffer if we redefine the estimates by assigning 
some default measure $Q(dx)$ in those rare cases when 
the search depth $\lambda(k(n),\l(n),J(n))$ exceeds the 
available record length $n$. It is difficult to say what 
the optimal growth path is for $k(n)$, $\l(n)$ and $J(n)$ 
without prior information about the spatial and temporal 
dependency structure of the process.

The special case of finite alphabet processes 
is most interesting and it is simpler because 
only 2 of the 3 parameters $k,\l,J$ play a role. 
We do not need an index for subfields of $\X$ because 
the obvious choice for $\B_k$ is the field of all 
subsets of $\X$. Also, it is convenient to choose 
the block length $\l_k$ equal to $k$ so that 
$\tau^k_j$ is the time for $j$ recurrences of $X^{-k}$.

In Section A we recall the ergodic theorem for
recurrence times that was derived by Wyner and 
Ziv \cite{bib36} and by Ornstein and Weiss \cite{bib24} for finite 
alphabet processes. In Section~B we define 
conditional probability mass function estimates 
$\hat P(x|X^{-n})$ and we prove consistency in 
mean if the block length $k(n)$ and the sample 
size $J_{k(n)}$ grow deterministically and 
sufficiently slowly with $n$. In Section~C we discuss generalizations 
for real-valued processes. 

\bigskip
\noindent
{\bf A. Recurrence Times}

\medskip
Let $\{X_t\}$ be a stationary ergodic process with 
values in a finite set $\X$. Starting at time $\tau^k_0=0$, 
the successive recurrence times $\tau^k_j$ of 
the $k$-block $X^{-k}$ are defined as follows:
\be\tau^k_j=\inf\{t>\tau^k_{j-1}:\,
(X_{-k-t},\ldots,X_{-1-t})=(X_{-k},\ldots,X_{-1})\}.\label{eq27}\ee
If $P\{X^{-k}=x^{-k}\}>0$ then by the results of 
Kac \cite{bib17} (see also Willems \cite{bib35}, 
Wyner and Ziv \cite{bib36}), 
\be E\{\tau^k_1|X^{-k}=x^{-k}\}
={1\over P\{X^{-k}=x^{-k}\}}.\label{eq28}\ee
Let $H$ denote the entropy rate of the stationary 
ergodic process $\{X_t\}$ in bits per symbol:
\be H=\lim_k\,-{1\over k}E\{\log P(X^k)\}
=\lim_k\,-{1\over k}E\{\log P(X^{-k})\}.\label{eq29}\ee
Wyner and Ziv \cite{bib36}, Theorem 3, invoked Kac's 
result and the Shannon-McMillan-Breiman theorem to prove 
that $\tau^k_1$ cannot grow faster than exponentially 
with limiting rate $H$ 
($\limsup_k\,k^{-1}\log\tau^k_1\le H$ almost surely).
Ornstein and Weiss \cite{bib24} then argued that $\tau^k_1$ 
will grow exponentially fast almost surely with 
limiting rate exactly equal to $H$:
\be k^{-1}\log\tau^k_1\to H\quad\hbox{almost surely.}\label{eq30}\ee 

Now suppose a sample of size $J_k$ is desired. The 
total time needed to find $J_k=J(k)\ge 1$ instances 
of the pattern $X^{-k}$ is equal to the recurrence 
time $\tau^k_{J(k)}$. The ratio $\tau^k_{J(k)}/J_k$ can 
be interpreted as the average inter-recurrence time: 
\be{\tau^k_{J(k)}\over J_k}={1\over J_k}
\sum_{1\le j\le J_k}(\tau^k_j-\tau^k_{j-1}).\label{eq31}\ee
We claim that like $\tau^k_j$, the average inter-recurrence 
time $\tau^k_{J(k)}/J_k$ cannot grow faster than 
exponentially with limiting rate $H$. The proof is 
based on Kac's result and the lemma that was developed 
by Algoet and Cover \cite{bib3} to give a simple proof of the 
Shannon-McMillan-Breiman theorem and a more general 
ergodic theorem for the maximum exponential growth rate 
of compounded capital invested in a stationary market.

\smallskip
\noindent
{\bf Theorem 2.}
{\sl Let $\{X_t\}$ be a stationary ergodic process with 
values in a finite set $\X$ and with entropy rate $H$ 
bits per symbol. If $\Delta_k=\Delta(k)$ is a sequence 
of numbers such that $\sum_k\,2^{-\Delta(k)}<\infty$,
then for arbitrary $J(k)=J_k>0$ we have
\be\log\left({\tau^k_{J(k)}\over J_k}\right) 
\le -\log P(X^{-k})+\Delta_k
\quad\hbox{eventually for large $k$,}\label{eq32}\ee
and consequently
\be\limsup_k\,{1\over k}\log\left({\tau^k_{J(k)}\over J_k}\right)
\le H\quad\hbox{almost surely.}\label{eq33}\ee
}

\noindent
{\sl Proof:}
The inter-recurrence times $\tau^k_j-\tau^k_{j-1}$
are identically distributed with the same conditional
distribution given $X^{-k}$ as the first recurrence 
time $\tau^k_1$.  By Kac's result, 
\be E\{\tau^k_{J(k)}|X^{-k}\}P(X^{-k})=J_k\,E\{\tau^k_1|X^{-k}\} P(X^{-k})
={J_k}.\label{eq34}\ee
(A referee pointed out that a result like this was 
also proved by Gavish and Lempel \cite{bib13}.) Thus the random 
variable $Z_k=P(X^{-k})\,\tau^k_{J(k)}/J_k$ has expectation
\be E\{Z_k\}=
E\left\{P(X^{-k})\,E\left\{{\tau^k_{J(k)}\over J_k}
\bigg|X^{-k}\right\}\right\}=1.\label{eq35}\ee
By the Markov inequality,
\be P\{\log Z_k>\Delta_k\}=P\{Z_k>2^{\Delta_k}\}
\le 2^{-\Delta_k}E\{Z_k\}=2^{-\Delta_k},\label{eq36}\ee
and by the Borel-Cantelli lemma $\log Z_k\le\Delta_k$
eventually for larger $k$. This proves (\ref{eq32}). 
Assertion (\ref{eq33}) follows from (\ref{eq32}) upon dividing both
sides by $k$ and taking the $\limsup$ as $k\to\infty$.
Indeed, $-k^{-1}\log P(X^{-k})\to H$ almost surely by the 
Shannon-McMillan-Breiman theorem and one may choose
$\Delta_k=2\log k$ so that $\Delta_k/k\to 0$.
\qed

It is worthwhile to observe that Theorem 2 can be 
generalized if the process $\{X_t\}$ is stationary 
but not necessarily ergodic. Let $P$ be a stationary 
distribution and let $P_\omega$ denote the ergodic 
mode of the actual process realization $\omega$. Then
by the ergodic decomposition theorem (see Theorem 7.4.1 
of Gray \cite{bib14}) and the monotone convergence theorem, 
\begin{eqnarray}
P\{X^{-k}=x^{-k}\} E\{\tau^k_{J(k)}|X^{-k}=x^{-k}\}&=&
\sum_{1\le t < \infty} t P\{X^{-k}=x^{-k}, \tau^k_{J(k)}=t\} \nonumber \\
&=&\sum_{1\le t < \infty} \int  t P_{\omega} \{X^{-k}=x^{-k}, 
\tau^k_{J(k)}=t\}P(d\omega) \nonumber\\
&=&\int \sum_{1\le t< \infty} t P_{\omega} \{X^{-k}=x^{-k}, 
\tau^k_{J(k)}=t\}P(d\omega) \nonumber \\
&=& \int  P_{\omega} \{X^{-k}=x^{-k}\} E_{\omega} 
\{\tau^k_{J(k)}|X^{-k}=x^{-k}\}P(d\omega) \nonumber \\
&=&\int J_k P(d\omega) \nonumber \\
&=& J_k.
\label{eq37}
\end{eqnarray}
It follows that $E\{P(X^{-k})\,\tau^k_{J(k)}\}=J_k$ and
\be\log(\tau^k_{J(k)}/J_k)\le -\log P(X^{-k})+\Delta_k
\quad\hbox{eventually for large $k$.}\label{eq38}\ee
The Shannon-McMillan-Breiman theorem for stationary
nonergodic processes asserts that $P(X^{-k})$ 
decreases exponentially fast with limiting rate 
$H(P_{\omega})$, so one may conclude that
\be\limsup_k\,{1\over k}\log\left({\tau^k_{J(k)}\over 
J_k}\right) \le H(P_{\omega}) \quad\hbox{almost surely.}\label{eq39}\ee
Thus the average inter-recurrence time $\tau^k_{J(k)}/J_k$ cannot 
grow faster than exponentially with limiting rate $H(P_{\omega})$, 
the entropy rate of the ergodic mode $P_{\omega}$.

\bigskip
\noindent
{\bf B. Conditional Probability Mass Function Estimates}

\medskip
In the finite alphabet case, the general estimator 
$\hat P_k(dx|X^{-\lambda(k)})$ that was defined 
in (\ref{eq16}) reduces to the conditional probability 
mass function estimate
\be\hat P_k(x|X^{-\lambda(k)})={1\over J_k}
\sum_{1\le j\le J_k}1\{X_{-\tau(k,j)}=x\}.\label{eq40}\ee
Here $k=\l_k$ is the block length and the sample 
size $J_k$ is monotonically increasing. The 
recurrence times $\tau^k_j$ of the $k$-block $X^{-k}$ 
were defined inductively for $j=1,2,3,\ldots$ in (\ref{eq27}).

We choose a slowly increasing sequence of block 
lengths $k(n)$ and set $\hat P(x|X^{-n})$ equal 
to $\hat P_{k(n)}(x|X^{-\lambda(k(n))})$ if this 
estimate can be computed from the available data 
segment $X^{-n}$. Otherwise, if $\lambda_{k(n)}>n$,
we truncate the search and define $\hat P(x|X^{-n})$
as the default measure $Q(x)=1/|\X|$. Thus for $n\ge 0$,
we define
\be\hat P(x|X^{-n})=\cases{\hat P_{k(n)}(x|X^{-\lambda(k(n))})
&if $\lambda(k(n))\le n$,\cr Q(x)&otherwise.\cr}\label{eq41}\ee
If $k(n)$ grows sufficiently slowly then truncation 
is a rare event and $\hat P(x|X^{-n})$ coincides most 
of the time with the weakly consistent estimator 
$\hat P_{k(n)}(x|X^{-\lambda(k(n))})$. The question 
is how fast the block length $k(n)$ and the sample 
size $J_{k(n)}$ may grow to get consistent estimates. 
To answer this question, we use our results about 
recurrence times. 

The inter-recurrence times $\tau^k_j-\tau^k_{j-1}$ 
have the same conditional distribution and hence 
the same conditional expectation given $X^{-k}$ as 
the first recurrence time $\tau^k_1$. The expected
inter-recurrence time is bounded as follows:
\be E\left\{{\tau^k_{J(k)}\over J_k}\right\}=E\{\tau^k_1\}
=\sum_{x^{-k}:\,P\{X^{-k}=x^{-k}\}>0}\,P\{X^{-k}=x^{-k}\}\,
E\{\tau^k_1|X^{-k}=x^{-k}\}\le |\X|^k.\label{eq42}\ee
If $\e_k>0$ then by the Markov inequality
\be P\left\{{\tau^k_{J(k)}\over J_k}>{|\X|^k\over\e_k}\right\}
\le \e_k.\label{eq43}\ee
If $\e_k\to 0$ then $P\{\tau^k_{J(k)}>J_k|\X|^k/\e_k\}\to 0$ 
and if $\sum_k\e_k<\infty$ then $\tau^k_{J(k)}\le J_k|\X|^k/\e_k$ 
eventually for large $k$ by the Borel-Cantelli lemma. 
This is similar to (\ref{eq32}) with $\e_k=2^{-\Delta(k)}$.
Since $\lambda(k)=k+\tau^k_{J(k)}$, we see that 
\be P\{\lambda(k(n))\le n\}\to 1\quad\hbox{as $n\to\infty$}\label{eq44}\ee 
if $J_k$ and $k(n)$ are chosen so that for some $\e_k>0$ with $\e_k\to 0$,
\be k(n)+J_{k(n)}|\X|^{k(n)}/\e_{k(n)}\le n
\quad\hbox{eventually for large $n$.}\label{eq45}\ee 
It suffices that $k(n)=(1-\epsilon)\log_{|\X|}n$ 
for some $0<\epsilon<1$ and $J_k=o(|\X|^{k\e/(1-\e)})$ 
so that $J_{k(n)}=o(n^\epsilon)$. (Noninteger values 
are rounded down to the nearest integer, as usual.) 
We can be slightly more aggressive.

\smallskip
\noindent
{\bf Theorem 3.}
{\sl Let $\{X_t\}$ be a stationary process with 
values in a finite set $\X$ and choose $Q(x)=|\X|^{-1}$ 
as default measure in (\ref{eq41}). If the block length
$k(n)$ and the sample size $J_{k(n)}$ are monotonically
increasing to infinity and satisfy
\be J_{k(n)}\,|\X|^{k(n)}=\O(n),\label{eq46}\ee
then the estimates $\hat P(x|X^{-n})$ in (\ref{eq41}) are consistent in mean:
\be\hat P(x|X^{-n})\to P(x|X^-)\quad\hbox{in $L^1(P)$.}\label{eq47}\ee
In particular, the estimates $\hat P(x|X^{-n})$ 
are consistent in mean if the block length is 
$k(n)= (1-\e) \log_{|\X|} n$ and the sample size 
is $J_{k(n)}=n^{\e}$ for some $0<\e<1$. 
}

\noindent
{\sl Proof:}
If the entropy rate $H$ is strictly less than $\log|\X|$ 
and $R$ is any constant such that $H<R<\log|\X|$ then 
by (\ref{eq33}), $\tau^k_{J(k)}$ is asymptotically bounded 
by $J_k 2^{Rk}$. It follows that 
\be \tau^{k(n)}_{J(k(n))}\le J_{k(n)}2^{Rk(n)}\le
J_{k(n)}|\X|^{k(n)}2^{(R-\log|\X|)k(n)}=o(n).\label{eq48}\ee
It is necessary for (\ref{eq46}) that $k(n)<\log_{|\X|}n$
eventually for large $n$ since $J_{k(n)}\to\infty$ by
assumption. Thus 
$\lambda(k(n))=k(n)+\tau^{k(n)}_{J(k(n))}=o(n)$ and
$\lambda(k(n))$ is upper bounded by $n$ eventually 
for large $n$.
If $H=\log |\X|$ then there is no guarantee that 
we can collect $J_{k(n)}$ samples from $X^{-n}$,
but the estimate $\hat P(x|X^{-n})$ will nevertheless
be consistent in mean if the default measure is
$Q(x)=|\X|^{-1}$ because the outcomes $X_t$ happen
to be independent identically distributed according 
to this distribution $Q(x)$ when $H=\log|\X|$.
\qed

The estimates $\hat P_k(x|X^{-\lambda(k)})$ in (\ref{eq40})
are consistent in the pointwise sense under certain conditions. 
For example, if $\{X_t\}$ is a stationary finite-state 
Markov chain with order $K$ then the empirical 
estimates $\hat P_k(x|X^{-\lambda(k)})$ are averages 
of bounded random variables $1\{X_{-\tau(k,j)}=x\}$
($j=1,2,\dots,J_k$) that are conditionally independent 
and identically distributed given $X^{-K}$ when $k\ge K$.  
It follows that the estimates $\hat P_k(x|X^{-\lambda(k)})$ 
converge exponentially fast in the number of samples $J_k$ 
to the conditional probability $P\{x|X^{-K}\}=P\{x|X^-\}$ 
and therefore the estimates are pointwise consistent. 
It is not known whether the estimates 
$\hat P_k(x|X^{-\lambda(k)})$ converge in the pointwise 
sense for all finite-alphabet stationary time series.

If we know the entropy rate $H$ in advance we can make use 
of it. In this case, weak consistency is guaranteed if
$k(n)=(1-\e)(\log n)/R$ for some $R>H$ and $J_{k(n)}=n^{\e}$. 
Indeed, if $H<r<R$ then $\lambda(k(n))<n$ eventually for 
large $n$ since
\begin{eqnarray} 
\lambda(k(n))&=& k(n)+\tau^{k(n)}_{J(k(n))}\nonumber\\
&\le& k(n)+J(k(n)) 2^{rk(n)}\nonumber\\
&=& \O(\log n)+n^{\e} n^{(1-\e)r/R}\nonumber\\
&=& o(n).
\end{eqnarray} 
If the entropy rate is not known in advance then
we must be prepared to deal with the worst case of 
nearly maximum entropy rate. The estimates will 
be wasteful if the entropy rate is low because they 
exploit only a small portion of the available data 
segment $X^{-n}$ when $H<\log|\X|$. If 
$k(n)=(1-\e)\log_{|\X|}n$ and $J_{k(n)}=n^\e$ 
then the length of the useful portion is about 
\be\tau^{k(n)}_{J(k(n))}\approx J_{k(n)}2^{Hk(n)} 
=n^{\e+(1-\e)H/\log|\X|}=n^\alpha,\label{eq49}\ee
where $\alpha=\e+(1-\e)H/\log|\X|$ varies linearly 
between $\e<\alpha\le 1$ as $0<H\le\log|\X|$.

The length $\lambda(k)=k+\tau^k_{J(k)}$ of the data 
record $X^{-\lambda(k)}$ that must be examined to 
collect $J_k$ samples of the pattern $X^{-k}$ grows 
approximately like $J_k 2^{Hk}$, which is polynomial 
in $J_k$ if $J_k$ grows exponentially fast with $k$. 
Also, the length $n$ of the segment $X^{-n}$ is just 
polynomial in the sample size $J_{k(n)}$ if $J_{k(n)}=n^\e$.
The strongly consistent estimates of Morvai, Yakowitz 
and Gy\"orfi \cite{bib21} are much less efficient: they collect 
$J$ samples from a data record whose length grows like 
a tower of exponentials in (\ref{eq7}). Their samples 
are very sparse because extremely stringent demands are 
placed on the context where those samples are taken.  
For the weakly consistent estimates of the present study, 
the demands on context are much less severe and so the 
samples are much more abundant although perhaps less 
trustworthy. Thus universal prediction is not hopelessly out 
of computational reach as it might seem for an algorithm 
whose input demands grow as a  tower of exponentials in (\ref{eq7}).

\bigskip
\noindent
{\bf C. Weak Consistency for Real-valued Processes}

\medskip
When $\X$ is the real line or a $\sigma$-compact Polish space,
the estimate $\hat P_k(dx|X^{-\lambda(k)})$ is defined by the 
formula in (\ref{eq16}). We now choose a nondecreasing unbounded sequence 
$k(n)$ and we define $\hat P(dx|X^{-n})$ as the empirical 
conditional distribution $\hat P_{k(n)}(dx|X^{-\lambda(k(n))})$ 
if this estimate can be computed from the available data 
segment $X^{-n}$. Otherwise, if $\lambda_{k(n)}>n$, 
we truncate the search and define $\hat P(dx|X^{-n})$ as 
some default measure $Q(dx)$. Thus
\be\hat P(dx|X^{-n})=\cases{\hat P_{k(n)}(dx|X^{-\lambda_{k(n)}})
&if $\lambda_{k(n)}\le n$,\cr Q(dx)&otherwise.\cr}\label{eq50}\ee
If $k(n)$ grows slowly then truncation is rare and 
$\hat P(dx|X^{-n})$ coincides most of the time with 
the estimator $\hat P_{k(n)}(dx|X^{-\lambda(k(n))})$ 
which is weakly consistent. The question is how slowly 
the partition index $k(n)$, the block length $\l(k(n))$
and the sample size $J(k(n))$ must grow with $n$ 
to get consistent estimates of $P(dx|X^-)$.
It suffices that $P\{\lambda(k(n))<n\}\to 1$.

\smallskip
\noindent
{\bf Theorem 4.} 
{\sl Let $\{X_t\}$ be a real-valued stationary ergodic 
time series and choose $\B_k$, $\l_k$ and $J_k$ as before. 
Let $\Xi_k$ denote the set of atoms of the finite field 
$\B_k$ and choose a nondecreasing unbounded sequence of
integers $k(n)$ and numbers $\e_k\to 0$ such that
\be n\ge \l_{k(n)}+J_{k(n)}|\Xi_{k(n)}|^{\l_{k(n)}}/\e_{k(n)}
\quad\hbox{eventually for large $n$.}\label{eq51}\ee
Then $P\{n\ge \lambda_{k(n)}\}\to 1$ as $n\to\infty$, and
the estimates $\hat P(dx|X^{-n})$ are weakly consistent:
for every set $B$ in the generating field $\bigcup_k\B_k$ 
we have
\be\hat P\{X\in B|X^{-n}\}\to P\{X\in B|X^-\}
\quad\hbox{in $L^1(P)$,}\label{eq52}\ee
and for every bounded continuous function $h(x)$ we have
\be\int h(x)\,\hat P(dx|X^{-n}) \to 
\int h(x)\,P(dx|X^-)\quad\hbox{in $L^1(P)$.}\label{eq53}\ee
}

\noindent
{\sl Proof:}
The inter-recurrence times $\tau^k_j-\tau^k_{j-1}$ 
($j=1,2,3,\ldots$) are identically distributed
conditionally given the pattern $[X^{-\l(k)}]^k$. 
By Kac's result, 
\be E\{\tau^k_{J(k)}|[X^{-\l_k}]^k\}
=J_k\,E\{\tau^k_1|[X^{-\l_k}]^k\}
={J_k\over P([X^{-\l_k}]^k)}.\label{eq54}\ee
It follows that
\be E\{\tau^k_{J(k)}\}=J_kE\{\tau^k_1\}
=\sum_{[x^{-\l(k)}]^k} P([x^{-\l(k)}]^k)
E\{\tau^k_{J(k)}|[x^{-\l(k)}]^k\}\le 
J_k|\Xi_k|^{\l(k)}.\label{eq55}\ee
(The sum is taken over $[x^{-\l(k)}]^k$ such that 
$P([x^{-\l(k)}]^k)=P\{[X^{-\l(k)}]^k=[x^{-\l(k)}]^k\}$
is strictly positive.) By the Markov inequality, 
\be P\{\lambda_k>\l_k+J_k|\Xi_k|^{\l(k)}/\e_k\}
=P\{\tau^k_{J(k)}>J_k|\Xi_k|^{\l(k)}/\e_k\}
\le {E\{\tau^k_{J(k)}\}\over J_k|\Xi_k|^{\l(k)}/\e_k}\le \e_k.\label{eq56}\ee
Assertions (\ref{eq52}) and (\ref{eq53}) follow from 
Theorem 1A and 1B because 
$P\{\l_k+J_k|\Xi_k|^{\l_k}/\e_k\ge \lambda_k\}\to 1$ 
and hence, in view of assumption (\ref{eq51}),
\be P\{n\ge \l_{k(n)}+J_{k(n)} {|\Xi_{k(n)}|}^{\l_{k(n)}}/\e_{k(n)}
\ge \lambda_{k(n)}\}\to 1\quad\hbox{as $n\to\infty$.}\label{eq57}\ee
This completes the proof of the theorem. \qed

The theorem remains valid in the stationary non-ergodic 
case. Indeed, let $P$ be a stationary distribution and 
let $P_{\omega}$ denote the ergodic mode of $\omega$. 
Then one may argue as above that 
$P_{\omega}\{\lambda_k\le\l_k+J_k|\Xi_k|^{\l(k)}/\e_k\}\to 1$. 
By the ergodic decomposition theorem
and Lebesgue's dominated convergence theorem,
\begin{eqnarray}
\lim_{k}P\{\lambda_k\le\l_k+J_k|\Xi_k|^{\l(k)}/\e_k\}
&=& \lim_{k}\int P_{\omega}
\{\lambda_k\le\l_k+J_k|\Xi_k|^{\l(k)}/\e_k\}P(d\omega)
\nonumber\\
&=& \int \lim_{k} P_{\omega}
\{\lambda_k\le\l_k+J_k|\Xi_k|^{\l(k)}/\e_k\}P(d\omega)
\nonumber\\
&=& \int 1 P(d\omega)=1.
\end{eqnarray}
Thus the conclusions of the theorem also hold 
for stationary nonergodic processes.

\bigskip
\centerline
{\large \bf IV. The Information Theoretic Point of View}

\medskip
In this section we discuss conditional distribution 
estimates $\hat P(dx|X^{-n})$ that are consistent in 
expected information divergence. Such estimates are 
also weakly consistent, but the converse is not 
necessarily true. It is possible to construct estimator 
sequences that are consistent in expected information 
divergence for all stationary processes with values 
in a finite alphabet, but not for all stationary 
processes with values in a countable infinite alphabet.   
There are connections with universal gambling or 
modeling schemes and with universal noiseless data 
compression algorithms for finite alphabet processes.
For more information on these subjects see Rissanen 
and Langdon \cite{bib28} and Algoet \cite{bib1}.

\medskip
\noindent
{\bf A. Consistency in Expected Information Divergence}

\smallskip
The Kullback-Leibler information divergence between 
two probability distributions $P$ and $Q$ on 
a measurable space $\X$ is defined as follows: 
if $P$ is dominated by $Q$ then
\be I(P|Q)=E_P\left\{\log\left({dP\over dQ}\right)\right\},\label{eq58}\ee
otherwise $I(P|Q)=\infty$. 
The variational distance is defined as
\be\|P-Q\|=\sup_{-1\le h(x)\le 1}\bigg|\int h\,dP-\int h\,dQ\bigg|,
\label{eq59}\ee
where the supremum is taken over all measurable 
functions $h(x)$ such that $|h(x)|\le 1$. 
If $p=dP/d\mu$ and $q=dQ/d\mu$ are the densities of
$P$ and $Q$ relative to a dominating $\sigma$-finite 
measure $\mu$ then $\|P-Q\|=\int |p-q|\,d\mu$.
Exercise 17 on p. 58 of Csisz\'ar 
and K\"orner \cite{bib11} asserts that
\be{\log e\over 2}\,\|P-Q\|^2\le I(P|Q).\label{eq61}\ee
It follows that $I(P|Q)\ge 0$ with equality iff $P=Q$.
Pinsker \cite{bib26}, pp. 13--15 proved the existence 
of a universal constant $\Gamma>0$ such that 
\be I(P|Q)\le E_P\left\{\bigg|\log\left({dP\over dQ}
\right)\bigg|\right\}\le I(P|Q)+\Gamma\sqrt{I(P|Q)}.\label{eq62}\ee
Barron \cite{bib6} simplified Pinsker's argument and proved 
that the constant $\Gamma=\sqrt{2}$ is best possible when 
natural logarithms are used in the definition of $I(P|Q)$.

Let $\{X_t\}$ be a stationary process with values 
in a complete separable metric space $\X$. The 
divergence between the true conditional distribution 
$P(dx|X^-)$ and an estimate $\hat P(dx|X^{-t})$ is 
a nonnegative function of the past $X^-$ which 
vanishes iff $P(dx|X^-)=\hat P(dx|X^{-t})$ 
$P$-almost surely. 
We say that the estimates $\hat P(dx|X^{-t})$ are 
consistent in information divergence for a class $\Pi$ 
of stationary distributions on $\X^\Z$ if for any $P\in\Pi$,
\be I(P_{X|X^-}|\hat P_{X|X^{-t}})\to 0
\quad\hbox{$P$-almost surely.}\label{eq63}\ee
We say that $\hat P(dx|X^{-t})$ is consistent in expected 
information divergence for the class $\Pi$ if for any 
$P\in\Pi$,
\be E_P\{I(P_{X|X^-}|\hat P_{X|X^{-t}})\}\to 0.\label{eq64}\ee
Such estimates are weakly consistent for all 
distributions in the class $\Pi$. Indeed, if $h(x)$ 
is any bounded measurable function on $\X$ 
with norm $\|h\|_\infty=\sup_x\,|h(x)|$ then
\be\big|\textstyle{\int} h(x)\,P(dx|X^-)
-\textstyle{\int}h(x)\,\hat P(dx|X^{-t})\big|
\le \|h\|_\infty\|P_{X|X^-}-\hat P_{X|X^{-t}}\|.\label{eq65}\ee
Applying the Csisz\'ar-Kemperman-Kullback inequality 
(\ref{eq61}), we see that
\be\bigg|\int h(x)\,P(dx|X^-)-\int h(x)\,\hat P(dx|X^{-t})
\bigg|^2\le {2\|h\|^2_\infty \over \log e}\,
I(P_{X|X^-}|\hat P_{X|X^{-t}}).\label{eq66}\ee
If $\hat P(dx|X^{-t})$ is consistent in expected information 
divergence for $\Pi$ then $\int h(x)\,\hat P(dx|X^{-t})$
converges in $L^2(P)$ and also in $L^1(P)$ to 
$\int h(x)\,P(dx|X^-)$ whenever $P\in\Pi$.

Suppose the outcomes $X_t$ are independent 
with identical distribution $P_X$ on $\X$. 
Barron, Gy\"orfi and van der Meulen \cite{bib7} have 
constructed estimates $\hat P(dx|X^{-t})$ that are 
consistent in information divergence and in expected 
information divergence when the true distribution $P_X$ 
has finite information divergence $I(P_X|M_X)<\infty$ 
relative to some known normalized reference measure $M_X$. 
Gy\"orfi, P\'ali and van der Meulen \cite{bib16} assume 
that $\X$ is the countable set of integers and argue that
for arbitrary conditional probability mass function 
estimates $\hat P(x|X^{-n})$, there exists some 
distribution $P_X$ with finite entropy such that
\be I(P_X|\hat P_{X|X^{-n}})=\infty\quad\hbox{almost surely for all $n$.}\label{eq67}\ee
Therefore, it is impossible to construct estimates 
$\hat P(dx|X^{-t})$ that are consistent in information
divergence or in expected information divergence for 
all independent identically distributed processes 
with values in an infinite space. For stationary 
processes with values in a finite alphabet, the 
constructions of Ornstein \cite{bib22} and Morvai, Yakowitz and Gy\"{o}rfi   
\cite{bib21} yield estimates $\hat P(x|X^{-t})$ such that 
$\log\hat P(x|X^{-t})$ converges almost surely to $\log P(x|X^-)$. 
It is still an open question as to whether these estimates 
are consistent in information divergence or whether 
modifications are needed to get such consistency. 
(The difficulty is that small changes in $\hat P(x|X^{-n})$
cause huge changes in $\log\hat P(x|X^{-n})$ when 
$\hat P(x|X^{-n})$ is small.) However, it is easy 
to construct estimates $\hat P(x|X^{-t})$ that are 
consistent in  expected information divergence.

\medskip
\noindent
{\bf B. Consistent Estimates for Finite-alphabet Processes}

\smallskip
Let $\{X_t\}$ be a stationary process with values 
in a finite set $\X$. We shall construct conditional
probability mass function estimates $\hat P(x|X^{-n})$ 
that are consistent in expected information divergence
for any stationary $P\in\P_s$.
Such estimates also converge to $P(x|X^-)$ in mean: 
for any stationary $P\in\P_s$ and $x\in\X$ we have
\be\hat P(x|X^{-n})\to P(x|X^-)\quad
\hbox{in $L^1(P)$.}\label{eq68}\ee
An observation of Perez \cite{bib25} implies that 
consistency in expected information divergence is 
equivalent to mean consistency of $\log \hat P(X|X^{-n})$.

\smallskip
\noindent
{\bf Theorem 5.}
{\sl Let $\{X_t\}$ be a stationary process with values 
in a finite alphabet $\X$. A sequence of conditional
probability mass function estimates $\hat P(x|X^{-n})$ 
is consistent in expected information divergence iff 
we have mean convergence
\be\log \hat P(X|X^{-n})\to \log P(X|X^-)
\quad\hbox{in $L^1$.}\label{eq69}\ee
}

\noindent
{\sl Proof:} 
Pinsker's inequality (\ref{eq62}) for $P(x|X^-)$ 
and $\hat P(x|X^{-n})$ asserts that
\begin{eqnarray}
I(P_{X|X^-}|\hat P_{X|X^{-n}})
&\le&
E\left\{\bigg|\log\left({P(X|X^-)\over \hat P(X|X^{-n})}\right)
\bigg|\bigg|X^-\right\} \nonumber\\
&\le& I(P_{X|X^-}|\hat P_{X|X^{-n}})
+\Gamma\sqrt{I(P_{X|X^-}|\hat P_{X|X^{-n}})}.
\end{eqnarray}
Taking expectations and using concavity 
of the square root function, we obtain
\begin{eqnarray}
E\{I(P_{X|X^-}|\hat P_{X|X^{-n}})\}
&\le& E\bigg|\log\left({P(X|X^-)\over \hat P(X|X^{-n})}\right)\bigg|
\nonumber\\
&\le& E\{I(P_{X|X^-}|\hat P_{X|X^{-n}})\}
+\Gamma\sqrt{E\{I(P_{X|X^-}|\hat P_{X|X^{-n}})\}}
\end{eqnarray}
by Jensen's inequality. 
This suffices to prove the theorem.
\qed

To construct the estimates $\hat P(x|X^{-n})$,
we start with probability mass functions $Q(x^n)$ 
on the product spaces $\X^n$ such that for every 
stationary distribution $P$ on $\X^\Z$,
\be n^{-1}I(P_{X^n}|Q_{X^n})\to 0\quad\hbox{as $n\to\infty$.}\label{eq73}\ee
Several methods are known for constructing such models 
$Q(x^n)$ -- see Section~C below. By Pinsker's inequality, 
convergence of the means in (\ref{eq73}) is equivalent 
to mean convergence
\be{1\over n}\log\left({P(X^n)\over Q(X^n)}\right)\to 0
\quad\hbox{in $L^1(P)$.}\label{eq74}\ee
Let now $Q(x|x^{-t})$ denote a shifted copy of 
the conditional probability mass function $Q(x_t|x^t)$ 
that appears in the chain rule expansion
$Q(x^n)=\prod_{0\le t<n}Q(x_t|x^t)$. The estimate 
$\hat P(x|X^{-n})$ is defined in terms of $Q(x^n)$ as 
\be\hat P(x|X^{-n})={1\over n}\sum_{0\le t<n}Q(x|X^{-t}).\label{eq75}\ee

\smallskip
\noindent
{\bf Theorem 6.}
{\sl Let $\X$ be a finite alphabet and let $\{Q(x^n)\}_{n\ge 1}$
be a model sequence such that (\ref{eq73}) or (\ref{eq74}) 
holds for all $P\in\P_s$. 
Then the conditional probability mass function estimates 
$\hat P(x|X^{-n})$ are consistent in expected information 
divergence for the class $\P_s$ of all stationary process 
distributions on $\X^\Z$.}

\noindent
{\sl Proof:}
The Kullback-Leibler divergence functional is 
convex in both arguments.  By the definition (\ref{eq75}) 
of $\hat P(x|X^{-n})$ and by Jensen's inequality, 
\be I(P_{X|X^-}|\hat P_{X|X^{-n}})\le {1\over n}
\sum_{0\le t<n}I(P_{X|X^-}|Q_{X|X^{-t}}).\label{eq76}\ee
Now we take expectations with respect to some 
distribution $P\in\P_s$. By stationarity and 
the chain rule expansion of information divergence, 
we obtain
\begin{eqnarray}
E_P\{I(P_{X|X^-}|\hat P_{X|X^{-n}})\}
&\le& {1\over n}\sum_{0\le t<n}E_P\{I(P_{X|X^-}|Q_{X|X^{-t}})\} \nonumber \\
&=& {1\over n}\sum_{0\le t<n}E_P\{I(P_{X_t|X^-X^t}|Q_{X_t|X^t})\} \nonumber \\
&=& {1\over n} E_P\{I(P_{X^n|X^-}|Q_{X^n})\} \nonumber \\
&=&{1\over n}E_P\{I(P_{X^n|X^-}|P_{X^n})\}
+{1\over n}I(P_{X^n}|Q_{X^n}).
\label{eq77}
\end{eqnarray}
Observe that
\be E_P\{I(P_{X^n|X^-}|P_{X^n})\}=H(X^n)-H(X^n|X^-)\label{eq78}\ee
where $H(X^n)=E_P\{-\log P(X^n)\}$ and 
$H(X^n|X^-)=E_P\{-\log P(X^n|X^-)\}$.
The entropy rate of the process is defined as
$H=H(X|X^-)=n^{-1}H(X^n|X^-)=\;\downarrow\lim_n\,n^{-1}H(X^n)$,
so one may conclude that
\be {1\over n}E_P\{I(P_{X^n|X^-}|P_{X^n})\}
={1\over n}H(X^n)-H \to 0\quad\hbox{as $n\to\infty$.}
\label{eq79}\ee
It follows from (\ref{eq77}) and (\ref{eq79})
that the estimates $\hat P(x|X^{-n})$ are consistent 
in expected information divergence, as claimed. 
\qed

The procedure which constructs $\hat P(x|x^{-n})$ 
from the models $Q(x^n)$ can be reversed. Indeed, 
let $\{\hat P(x|x^{-t}\}_{t\ge 0}$ be a sequence
such that for every stationary distribution $P\in\P_s$,
the expected information divergence of $P(x|X^-)$ 
relative to $\hat P(x|X^{-t})$ is finite for all 
$t$ and vanishes in the limit as $t\to\infty$.
Let $\hat P(x_t|x^t)$ be constructed from the 
$t$-past at time $t$ in the same way as $\hat P(x|x^{-t})$ 
was constructed from the $t$-past at time $0$. 
The Kullback-Leibler information divergence of the true 
marginal distribution $P(x^n)$ with respect to the compounded 
model $\hat P(x^n)=\prod_{0\le t<n}\hat P(x_t|x^t)$
admits the chain rule expansion
\be I(P_{X^n}|\hat P_{X^n})=\sum_{0\le t<n}
E_P\{I(P_{X_t|X^t}|\hat P_{X_t|X^t})\}.\label{eq93}\ee
By stationarity 
\be E_P\{I(P_{X_t|X^t}|\hat P_{X_t|X^t})\}
=E_P\{I(P_{X|X^{-t}}|\hat P_{X|X^{-t}})\}.\label{eq94}\ee
The divergence of $P(x|X^{-t})$ relative to $\hat P(x|X^{-t})$
is bounded by the divergence of $P(x|X^-)$ relative 
to $\hat P(x|X^{-t})$ since we have the decomposition
\be I(P_{X|X^-}|\hat P_{X|X^{-t}})=I(P_{X|X^-}|P_{X|X^{-t}})
+I(P_{X|X^{-t}}|\hat P_{X|X^{-t}}).\label{eq95}\ee 
From (\ref{eq93}), (\ref{eq94}) and (\ref{eq95}) 
one may conclude that
\be I(P_{X^n}|\hat P_{X^n})\le 
\sum_{0\le t<n}E_P\{I(P_{X|X^-}|\hat P_{X|X^{-t}})\}.\label{eq96}\ee
If the expected divergence between $P(x|X^-)$ and 
$\hat P(x|X^{-t})$ is finite and vanishes in the 
limit as $t\to\infty$ then the models 
$\hat P(x^n)=\prod_{0\le t<n}\hat P(x_t|x^t)$
have vanishing expected per-symbol divergence:
for all $P\in\P_s$ we have
\be n^{-1}I(P_{X^n}|\hat P_{X^n})\to 0.\label{eq97}\ee

The results of Shields \cite{bib34} imply that
there can be no universal bound on the speed of 
convergence in expected information divergence.
Indeed, if the expected divergence of $P(x|X^-)$ relative 
to $\hat P(x|X^{-t})$ were always $\O(\beta_t)$ where 
$\beta_t\to 0$, then we could construct a modeling 
scheme $\{\hat P(x_t|x^t)\}_{t\ge 0}$ such that
the divergence of $P(x^n)$ relative to 
$\hat P(x^n)=\prod_{0\le t<n}\hat P(x_t|x^t)$ would 
be $\O(\beta_0+\ldots+\beta_{n-1})$. The per-symbol 
divergence of $P(x^n)$ relative to $\hat P(x^n)$ 
would vanish with universal rate 
$\O[n^{-1}(\beta_0+\ldots+\beta_{n-1})]$, which is
impossible.

To obtain bounds on the per-symbol divergence 
one must restrict the process distribution to some 
manageable class. In particular, suppose $\Pi$ is a class 
of Markov processes that is smoothly parametrized by $k$ 
free parameters and consider models $\hat P(x^n)$ for 
which the per-symbol divergence attains Rissanen's 
\cite{bib27} lower bound: 
\be{1\over n}I(P_{X^n}|\hat P_{X^n})=
{k\log n\over 2n}(1+o(1)).\label{eq98}\ee
If we set $Q(x^n)=\hat P(x^n)$ and define 
$\hat P(x|X^{-n})$ as in (\ref{eq75}), then (\ref{eq77})  reduces to the bound  
\be E_P\{I(P_{X|X^-}|\hat P_{X|X^{-n}})\}\le {k\log n\over 2n}
(1+o(1)), \quad P\in\Pi.\label{eq99}\ee
It is often possible to construct a prequential modeling scheme 
$\{\hat P(x_t|x^t)\}_{t\ge 0}$ such that the expected 
divergence of $P(x_t|X^t)$ relative to $\hat P(x_t|X^t)$ 
vanishes like $(k\log e)/(2t)$ for all process distributions in the class  $\Pi$. 
An incremental bound of order $(k\log e)/(2t)$ yields a normalized cumulative bound of order 
$n^{-1} \sum_{t<n} (k\log e)/(2t) \approx(k\log n)/(2n)$. 
By shifting $\hat P(x_n|X^n)$ we obtain estimates 
$\hat P(x|X^{-n})$ such that the expected divergence of $P(x|X^-)$ relative to $\hat P(x|X^{-n})$ 
vanishes like $(k\log e)/(2n)$. This bound of order 
$(k \log e)/(2n)$ for $\hat P(x|X^{-n})$ is clearly better than the bound $(k\log n)/(2n)$.

\bigskip
\noindent
{\bf C. Modeling and Data Compression}

\medskip
Any universal data compression scheme for stationary 
processes with finite alphabet $\X$ can be used as a basis 
for the construction of models $Q(x^n)$ satisfying (\ref{eq73}) 
or (\ref{eq74}). Indeed, let $l(x^n)$ denote the length of a 
uniquely decipherable block-to-variable-length binary 
code for sequences $x^n\in\X^n$. The redundancy of the 
code for $X^n$ is defined as the difference between 
the actual codeword length $l(X^n)$ and the ideal 
description length $-\log P(X^n)$:
\be r(X^n)=l(X^n)+\log P(X^n).\label{eq80}\ee
The expected redundancy $E_P\{r(X^n)\}$ is equal to the 
information divergence between the true probability 
mass function $P(x^n)$ and the model
\be Q'(x^n)=2^{-l(x^n)},\quad x^n\in\X^n.\label{eq81}\ee
For a universal noiseless coding scheme, the expected 
per-symbol redundancy will vanish: 
\be{1\over n}E_P\{r(X^n)\}={1\over n}I(P_{X^n}|Q'_{X^n})
\to 0\quad\hbox{for all $P\in\P_s$.}\label{eq82}\ee
The model $Q'(x^n)$ is not necessarily normalized but is
always a subprobability measure, by the Kraft-McMillan inequality.
However, if (\ref{eq82}) holds for a sequence of subnormalized 
models $Q'(x^n)$ then (\ref{eq82}) will certainly hold for the 
normalized models
\be Q(x^n)={Q'(x^n)\over\sum_{\xi^n\in\X^n}Q'(\xi^n)}.\label{eq83}\ee

Theorem 4 of Algoet \cite{bib1} implies that for any stationary 
ergodic distribution $P$, the per-symbol description 
length of uniquely decipherable codes is asymptotically 
bounded below almost surely by the entropy rate 
$H(P)=\lim_n\,n^{-1}E_P\{-\log P(X^n)\}$:  
\be\liminf_n\,n^{-1}l(X^n)\ge H(P)\quad\hbox{$P$-almost surely.}\label{eq85}\ee
It is well known that there exist universal noiseless 
codes for which the per-symbol description length 
almost surely approaches the entropy rate of the ergodic mode $P_{\omega}$ 
with probability one under any stationary distribution $P$:
\be n^{-1}l(X^n(\omega))\to H(P_{\omega})\quad\hbox{$P$-almost surely, for all $P\in\P_s$.}\label{eq86}\ee
This is true in particular for the data compression
algorithm of Ziv and Lempel \cite{bib38}, by Theorem 12.10.2 
of Cover and Thomas \cite{bib10} or by the results of Ornstein 
and Weiss \cite{bib24}. Other examples of noiseless codes
satisfying (\ref{eq86}) for every stationary ergodic $P$ 
have been proposed by Ryabco \cite{bib29}, Ornstein 
and Shields \cite{bib23}, and Algoet \cite{bib1}.
Choosing the best among the given code with length $l(x^n)$ 
and a fixed-length code with length $\lceil n\log|\X|\rceil$ 
and adding one bit of preamble to indicate which code is better, 
one obtains a uniquely decipherable code with length
\be l'(x^n)=1+\min\{l(x^n),\lceil n\log|\X|\rceil\}.\label{eq84}\ee
The codeword may expand by one bit, but the per-symbol 
description length is now bounded by $\log|\X|+2n^{-1}$ 
and (\ref{eq86}) holds universally not only in the 
pointwise sense but also in mean. The corresponding 
models $Q(x^n)$ are universal in the sense that for 
any $P\in\P_s$,
\be{1\over n}\log\left({P(X^n)\over Q(X^n)}\right)
\to 0\quad\hbox{$P$-almost surely and in $L^1(P)$.}\label{eq87}\ee

Ryabco \cite{bib29} and Algoet \cite{bib1} have 
constructed probability measures $Q$ with marginals 
$Q(x^n)$ such that the pointwise convergence in 
(\ref{eq87}) holds for every stationary $P\in\P_s$. 
Each marginal $Q(x^n)$ is equal to the compounded product 
$Q(x^n)=\prod_{0\le t<n}Q(x_t|x^t)$, and Ryabco's
scheme has the extra property that when $P$ is 
finite order Markov, 
\be\log\left({P(X_t|X^t)\over Q(X_t|X^t)}\right)\to 0
\quad\hbox{$P$-almost surely.}\label{eq88}\ee
Rissanen and Langdon \cite{bib28} and Langdon \cite{bib18} 
previously observed that the Lempel-Ziv algorithm defines a 
sequential predictive modeling scheme $Q=\{Q(x_t|x^t)\}$. 
The per-symbol divergence vanishes pointwise in 
the Ces\`aro mean sense, for every $P\in\P_s$:
\be{1\over n}\log\left({P(X^n)\over Q(X^n)}\right)
={1\over n}\sum_{0\le t<n}
\log\left({P(X_t|X^t)\over Q(X_t|X^t)}\right)\to 0
\quad\hbox{$P$-almost surely}\label{eq91}.\ee
However, the pointwise convergence in (\ref{eq88}) 
must fail for some $P\in\P_s$ because the quality of 
the predictive model $Q(x_t|X^t)$ degrades whenever the 
Lempel-Ziv incremental parsing procedure comes to the end 
of a phrase. 
The leaves of the dictionary tree and the nodes with few 
descendants are exactly those where empirical evidence is 
still lacking to make a reliable forecast. The number of 
times a node has been visited is equal to the number of 
leaves in the subtree rooted at that node, and if this 
number is small then the predictive model for the next 
symbol is a poor estimate based on few samples.

If the estimates $\hat P(x|X^{-t})$ are universally 
consistent in expected information divergence then 
$\log[P(X|X^-)/\hat P(X|X^{-t})]\to 0$ in $L^1(P)$ for
all stationary $P\in\P_s$ by Theorem~5. Thus the shifted
estimates $\hat P(x_t|X^t)$ are universally consistent 
in the sense that for all $P\in\P_s$,
\be\log\left({P(X_t|X^t)\over \hat P(X_t|X^t)}\right)\to 0
\quad\hbox{in $L^1(P)$.}\label{eq915}\ee
Bailey \cite{bib5} and Ryabco \cite{bib30} proved 
that no modeling scheme $Q$ exists such that the 
pointwise convergence in (\ref{eq88}) holds for 
every stationary ergodic distribution $P$. The argument 
of \cite{bib30} shows that for any modeling scheme 
$Q$ there exists a stationary ergodic distribution 
$P$ on $\X^\Z$ where $\X=\{a,b,c\}$ such that $P$ 
fails to satisfy both (\ref{eq88}) and the statement
\be P(X_t|X^t)-Q(X_t|X^t)\to 0\quad\hbox{$P$-almost surely.}\label{eq92}\ee
The offending $P$ is determined by a Markov chain with 
a countable set of states $\{0,1,2,\ldots\}$. Given that 
the Markov chain is in state $i$, it moves to state $0$ 
with probability $1/2$ and generates the letter $a$, 
or it moves to state $i+1$ with probability $1/2$ and 
generates the letter $b$ or $c$ with conditional 
probability $\Delta_i$ and $(1-\Delta_i)$, where 
$\Delta_i$ is a parameter equal to either $1/3$ or $2/3$. 
The distribution of the Markov chain is determined by 
the infinite sequence $\Delta=(\Delta_0,\Delta_1,\ldots)$.
If the Markov chain is started in its stationary 
distribution then the resulting distribution $P_\Delta$ 
on the sequence space $\X^\infty$ is stationary ergodic. 
Exact prediction is impossible when the Markov chain 
visits a state $i$ which it has not visited before, 
because the predictor doesn't know whether the probability 
$\Delta_i/2$ of next seeing symbol $b$ is equal to $1/3$ or $1/6$. 
The Markov chain will visit states with arbitrarily large 
labels $i$, and the predictor must make inaccurate predictions 
infinitely often with positive probability  under distribution 
$P_\Delta$ for some $\Delta=(\Delta_0,\Delta_1,\ldots)$.

\bigskip
\centerline{\large \bf V. Application to Online Prediction}

\medskip
In this section we discuss some applications of 
the estimates $\hat P(dx|X^{-t})$ to on-line
prediction, regression and classification. 
We deal with special cases of a sequential decision 
problem that can be formulated abstractly as follows.

Let $\{X_t\}$ be a stationary process with values in 
the space $\X$ and let $l(x,a)$ be a loss function on 
$\X\times\A$ where $\A$ is a space of possible actions. 
We assume that $\X$ is a complete and $\A$ is a compact 
separable metric space and the loss function $l(x,a)$ 
is bounded and continuous on $\X\times\A$.
We wish to select nonanticipating actions $A_t=A_t(X^t)$ 
with knowledge of the past $X^t=(X_0,\ldots,X_{t-1})$ so 
as to minimize the long run average loss per decision:
\be\limsup_n\,{1\over n}\sum_{0\le t<n}l(X_t,A_t)
={\rm Min!}\label{eq100}\ee

If the process distribution is known a priori 
then the optimum strategy is to select actions 
$A^*_t=\arg\min_{a\in\A}\,E\{l(X_t,a)|X^t\}$ that 
attain the minimum conditional expected loss given 
the available information $X^t$ at each time $t$. 
Suppose $P$ is stationary and let $L(X_t|X^t)$
denote the expectation of the minimum conditional 
expected loss given the $t$-past at time $t$:
\be L(X_t|X^t)=E\{l(X_t,A^*_t)\}
=\inf_{A_t=A_t(X^t)}\,E\{l(X_t,A_t)\}.\label{eq101}\ee
Similarly let $L(X|X^{-t})$ and $L(X|X^-)$ denote 
the minimum expected loss given the $t$-past and 
the minimum expected loss given the infinite past 
at time $0$. By stationarity $L(X_t|X^t)=L(X|X^{-t})$, 
and $L(X|X^{-t})$ is clearly monotonically decreasing 
to a limit which by continuity must be $L(X|X^-)$. 
Thus for any stationary distribution $P$ one may define
\be L^*(P)=\;\downarrow\lim_t\,L(X_t|X^t)
=\;\downarrow\lim_t\,L(X|X^{-t})=L(X|X^-).\label{eq102}\ee
If $P$ is stationary ergodic then the minimum long 
run average loss is well defined and almost surely 
equal to $L^*(P)=L(X|X^-)$ by Theorem~6 of Algoet \cite{bib2}:
\be{1\over n}\sum_{0\le t<n}l(X_t,A^*_t)\to L^*(P)
\quad\hbox{$P$-almost surely and in $L^1(P)$.}\label{eq103}\ee

Now suppose the process distribution is unknown a priori. 
It is shown in Section~V.B of Algoet \cite{bib2} that there exist 
nonanticipating actions $\hat A^*_t=\hat A^*_t(X^t)$ 
which attain the minimum long run average loss $L^*(P)$ 
with probability one under any stationary ergodic process 
distribution $P$ on $\X^\Z$. The actions $\hat A^*_t$ are 
constructed by a plug-in approach as follows. Choose 
estimates $\hat P(dx|X^{-t})$ that converge in law to the 
conditional distribution $P(dx|X^-)$ with probability one 
under any stationary $P$ and construct $\hat P(dx_t|X^t)$ 
from $X^t$ in the same way as $\hat P(dx|X^{-t})$ was 
computed from $X^{-t}$. Then $\hat A^*_t$ is defined as 
an action that attains the minimum conditional expected 
loss given $X^t$ under $\hat P(dx_t|X^t)$:
\be\hat A^*_t=\arg\min_{a\in\A}\,
\int l(x_t,a)\,\hat P(dx_t|X^t).\label{eq106}\ee
The average loss incurred by the actions $\hat A^*_t$ 
converges pointwise to the minimum long run average loss 
$L^*(P)$. 

In this paper we rely on conditional distribution 
estimates $\hat P(dx|X^{-t})$ that are weakly 
consistent but hopefully more efficient than the 
pointwise consistent estimates of \cite{bib22}, \cite{bib1}, \cite{bib21}.
We limit our attention to certain on-line prediction 
problems, when $\X=\A$ is a compact separable metric 
space and the loss $l(x,\hat x)$ is a continuous 
increasing function of the distance between the 
outcome $x$ and the prediction $\hat x$. 
In classification problems $\X=\A$ is a finite set, 
$l(x,\hat x)=1\{x\neq \hat x\}$ is the Hamming distance, 
and we wish to predict each outcome $X_t$ with knowledge 
of the past $X^t$ so as to minimize the long run average 
rate of incorrect guesses. In regression problems $\X$ 
is a finite closed interval, $l(x,\hat x)$ is the 
squared euclidean distance, and the goal is to predict 
$X_t$ from the past $X^t$ so that the long run average 
of the squared prediction error is smallest possible. 
We show that if the estimates $\hat P(dx|X^{-t})$ are 
weakly consistent, then the minimum long run average 
loss in regression and classification is universally 
attained in the sense of mean convergence in $L^1(P)$. 
The proof is based on the following generalization of 
von Neumann's mean ergodic theorem, which parallels 
Breiman's \cite{bib8} generalization of Birkhoff's pointwise
ergodic theorem. See also Perez \cite{bib25}.

\smallskip
\noindent
{\bf Lemma.}
{\sl Suppose $(\Omega,\F,P,T)$ is a stationary ergodic system.
If $g$ and $\{g_t\}_{t\ge 0}$ are integrable random variables such 
that $g_t\to g$ in $L^1(P)$, then 
\be{1\over n}\sum_{0\le t<n}g_t\circ T^t\to E\{g\}\quad
\hbox{in $L^1(P)$.}\label{eq107}\ee
}

\noindent
{\sl Proof:} The mean ergodic theorem asserts that
\be{1\over n}\sum_{0\le t<n}g\circ T^t\to E\{g\}\quad
\hbox{in $L^1(P)$,}\label{eq108}\ee
and it is clear that 
\be{1\over n}\sum_{0\le t<n}[g_t\circ T^t-g\circ T^t]
\to 0\quad \hbox{in $L^1(P)$}\label{eq109}\ee
since the triangle inequality, stationarity 
and the assumption $E|g_t-g|\to 0$ imply that 
\be E\bigg|{1\over n}\sum_{0\le t<n}[g_t\circ T^t-g\circ T^t] \bigg|
\le {1\over n}\sum_{0\le t<n}E|g_t\circ T^t-g\circ T^t|
={1\over n}\sum_{0\le t<n}E|g_t-g|\to 0.\label{eq110}\ee
Addition of (\ref{eq108}) and (\ref{eq109}) yields (\ref{eq107}).
\qed

\bigskip
\noindent
{\bf A. Regression}

\medskip
Let $\{X_t\}$ be a stationary ergodic real-valued time series 
with finite variance.  We wish to predict each outcome $X_t$ 
with knowledge of the past $X^t$ so that the squared prediction 
error $|X_t-\hat X_t|^2$ is smallest possible in the long run 
average sense.  The minimum long run average is equal to 
the minimum mean squared error given the infinite past, 
that is the variance of the innovation $X-E\{X|X^-\}$.
If the outcomes $X_t$ are independent and identically distributed 
then the sample mean $\hat X_t=(X_0+\ldots+X_{t-1})/t$ is 
an optimal estimator in the long run. It is challenging to 
construct on-line predictors $\hat X_t$ that asymptotically 
attain the minimum squared prediction error in a universal 
sense for all stationary ergodic real-valued processes with 
finite variance. Here, we consider the simple case of 
stationary processes with values in a finite interval 
$\X=[-K,K]$. We do not assume that $K$ is known a priori.

Let $\{\hat P(dx|X^{-t})\}_{t\ge 0}$ denote a weakly 
consistent sequence of conditional distribution estimates 
as in Section III. Since $h(x)=x$ is a bounded continuous 
function on $\X$, it follows from Theorem 4 that 
$\hat X_{-t}\to \hat X$ in probability where
\be\hat X_{-t}=\int x\,\hat P(dx|X^{-t}),\quad
\hat X=E\{X|X^-\}=\int x\,P(dx|X^-).\label{eq111}\ee
Note that $\hat X_{-t}$ is not an estimate of $X_{-t}$
but an estimate of $X=X_0$ based on the $t$-past $X^{-t}$.
At time $t$ we consider the conditional distribution
estimate $\hat P(dx_t|X^t)$ and the predictor
\be\hat X_t=\int x_t\,\hat P(dx_t|X^t).\label{eq112}\ee
By construction $\hat X_t$ is the sample mean of 
some subset of the past outcomes $X_0,\ldots,X_{t-1}$,
except in rare cases when $\hat X_t$ is equal to the
default value $\int x\,Q(dx)$. The obvious choice for 
$Q(dx)$ is the Dirac measure that places unit mass at 
$x=0$, so that $\int x\,Q(dx)=0$. For any stationary 
ergodic process distribution $P$ on $\X^\Z$ we have
\be|X-\hat X_{-t}|^2\to |X-\hat X|^2\quad\hbox{in $L^1(P)$,}\label{eq113}\ee
and consequently, by the Lemma,
\be{1\over n}\sum_{0\le t<n}|X_t-\hat X_t|^2
\to E|X-\hat X|^2\quad\hbox{in $L^1(P)$.}\label{eq114}\ee

\bigskip
\noindent
{\bf B. On-line Prediction and Classification}

\medskip
Let $\{X_t\}$ be a random process with values in 
a finite set $\X$. We wish to predict the outcomes 
$X_t$ with knowledge of the past $X^t$ so as to 
minimize the long run average rate of incorrect 
guesses. The best predictor for $X=X_0$ given 
the infinite past $X^-$ is given by
\be\hat X=\arg\max_{x\in\X}\,P\{X=x|X^-\}.\label{eq115}\ee
If the process distribution $P$ is stationary ergodic 
then the minimum long run average rate of prediction 
errors is equal to the error probability
\be P\{X\neq\hat X\}=1-E\{P\{X=\hat X|X^-\}\}
=1-E\left\{\max_{x\in\X}\,P\{X=x|X^-\}\right\}.\label{eq116}\ee
If the process distribution is unknown, we choose 
some conditional probability mass estimates 
$\hat P\{X=x|X^{-t}\}$ that converge in mean 
to $P\{X=x|X^-\}$ for every stationary process 
distribution $P\in\P_s$ and $x\in\X$:
\be\hat P\{X=x|X^{-t}\}\to P\{X=x|X^-\}
\quad\hbox{in $L^1(P)$.}\label{eq117}\ee
We construct $\hat P\{X_t=x|X^t\}$ from the past $X^t$ 
in the same way as $\hat P\{X=x|X^{-t}\}$ was computed
from $X^{-t}$ and we define the predictor
\be\hat X_t=\arg\max_{x\in\X}\,\hat P\{X_t=x|X^t\},\label{eq118}\ee

\smallskip
\noindent
{\bf Theorem 7.}
{\sl Let $\{X_t\}$ be a stationary ergodic process with 
values in a finite set $\X$. If the conditional
probability estimates $\hat P\{X=x|X^{-t}\}$ are 
weakly consistent, then the predictor $\hat X_t$ 
achieves the minimum long run average rate of incorrect 
guesses in probability. Thus for any stationary ergodic 
distribution $P$ on $\X^\Z$ we have mean convergence
\be{1\over n}\sum_{0\le t<n}1\{X_t\neq \hat X_t\}\to
P\{X\neq \hat X\}\quad\hbox{in $L^1(P)$.}\label{eq119}\ee
}

\noindent
{\sl Proof:}
Observe that $\hat X_t(\omega)=\hat X_{-t}(T^t\omega)$
where $T$ is the left shift on $\X^\Z$ and where
\be\hat X_{-t}=\arg\max_{x\in\X}\,\hat P\{X=x|X^{-t}\}.\label{eq120}\ee
For any stationary ergodic $P$ we have, by 
weak consistency of $\hat P\{X=x|X^{-t}\}$ 
and continuity of the maximum function, 
\be\max_{x\in\X}\,\hat P\{X=x|X^{-t}\}\to
\max_{x\in\X}\,P\{X=x|X^-\}\quad\hbox{in probability}\ \label{eq121}\ee
or equivalently
\be\hat P\{X=\hat X_{-t}|X^{-t}\}\to
P\{X=\hat X|X^-\}\quad\hbox{in $L^1(P)$.}\label{eq122}\ee
Since $[P\{X=x|X^-\}-\hat P\{X=x|X^{-t}\}]\to 0$ 
in $L^1(P)$ by weak consistency and
\be|P\{X=\hat X_{-t}|X^-\}-\hat P\{X=\hat X_{-t}|X^{-t}\}|
\le \sum_{x\in\X}|P\{X=x|X^-\}-\hat P\{X=x|X^{-t}\}|,\label{eq123}\ee
we see that
\be[P\{X=\hat X_{-t}|X^-\}-\hat P\{X=\hat X_{-t}|X^{-t}\}]
\to 0\quad\hbox{in $L^1(P)$.}\label{eq124}\ee
It follows from (\ref{eq122}) and (\ref{eq124}) that 
\be P\{X=\hat X_{-t}|X^-\}\to P\{X=\hat X|X^-\}
\quad\hbox{in $L^1(P)$}\label{eq125}\ee
and consequently, by the Lemma,
\be{1\over n}\sum_{0\le t<n}P\{X_t\neq\hat X_t|X^-X^t\}
\to E\{P\{X\neq \hat X|X^-\}\}=P\{X\neq\hat X\}
\quad\hbox{in $L^1(P)$.}\label{eq126}\ee
Now observe that 
\be\Delta_t=1\{X_t\neq \hat X_t\}-P\{X_t\neq\hat X_t|X^-X^t\}\label{eq127}\ee
is a bounded martingale difference sequence with 
respect to the $\sigma$-fields $\sigma(X^-X^t)$ and hence 
\be{1\over n}\sum_{0\le t<n}\Delta_t=
{1\over n}\sum_{0\le t<n}
[1\{X_t\neq \hat X_t\}-P\{X_t\neq\hat X_t|X^-X^t\}]
\to 0\quad\hbox{in $L^1(P)$}\label{eq128}\ee
(and also $P$-almost surely). In fact, the Ces\`aro
means of $\Delta_t$ vanish exponentially fast by 
Azuma's \cite{bib4} exponential inequalities for bounded
martingale differences.  Addition of (\ref{eq126}) and (\ref{eq128}) 
yields the conclusion (\ref{eq119}).
\qed

Feder, Merhav and Gutman \cite{bibFeder} used the Lempel-Ziv algorithm as a method 
for sequential prediction of individual sequences. 

\bigskip
\noindent
{\bf C. Problems with Side Information}

\medskip
A well studied problem in statistical decision theory,
pattern recognition and machine learning is to infer
the class label $X_t$ of an item at time $t$ from 
a covariate or feature vector $Y_t$ and a training 
set $X^tY^t=(X_0,Y_0,\ldots,X_{t-1},Y_{t-1})$. It 
is often reasonable to assume that the successive 
pairs $(X_t,Y_t)$ are independent and identically
distributed, but sometimes defective items tend 
to come in batches or in periodic runs and in those 
cases it may be profitable to exploit dependencies 
between new items and recent or not so recent items. 
Here we assume that the pair process $\{(X_t,Y_t)\}$ 
is stationary ergodic and we try to exploit statistical 
dependencies of arbitrarily long range, although we have 
no idea what kind of dependencies to expect a priori. 
The minimum long run average misclassification rate 
is again equal to $P\{X\neq\hat X\}$, but now $\hat X$ 
is the best predictor of $X=X_0$ given the infinite past
$X^-Y^-=(\ldots,X_{-2},Y_{-2},X_{-1},Y_{-1})$ and
the side information $Y=Y_0$:
\be\hat X=\arg\max_{x\in\X}\,P\{X=x|X^-Y^-Y\}.\label{eq129}\ee
The minimum misclassification rate will be asymptotically
attained in probability by the predictors
\be\hat X_t=\arg\max_{x\in\X}\,\hat P\{X_t=x|X^tY^tY_t\},\label{eq130}\ee
where $\hat P\{X_t=x|X^tY^tY_t\}$ is a shifted 
version of a conditional probability estimate 
$\hat P\{X=x|X^{-t}Y^{-t}Y\}$ such that for any 
stationary process distribution $P$ on $(\X\times\Y)^\Z$,
\be\hat P\{X=x|X^{-t}Y^{-t}Y\}\to P\{X=x|X^-Y^-Y\}
\quad\hbox{in $L^1(P)$.}\label{eq131}\ee
Such estimates $\hat P\{X=x|X^{-t}Y^{-t}Y\}$ can be 
constructed by generalizing the methods of Sections~II and III.

In fact, let $\X$ and $\Y$ be complete separable metric spaces 
and let $\{\B_k\}_{k\ge 1}$ and $\{\C_k\}_{k\ge 1}$
be increasing sequences of finite subfields that
asymptotically generate the Borel $\sigma$-fields on 
$\X$ and $\Y$. We assume that $\X$ is $\sigma$-compact
and the fields $\B_k$ are constructed as in the 
paragraph after Theorem 1B. Let $[x]^k$ and $[y]^k$ 
denote the atoms of $\B_k$ and $\C_k$ that contain 
the points $x\in\X$ and $y\in\Y$, and consider the
sequence of past recurrence times $\tau^k_j=\tau(k,j)$ 
of the pattern $[X^{-\l(k)}Y^{-\l(k)}Y]^k$. 
Then for every stationary process distribution $P$ 
on $(\X\times\Y)^\Z$, 
\be\hat P(dx|X^{-\lambda_k}Y^{-\lambda_k}Y)=
{1\over J_k}\sum_{1\le j\le J_k}\delta_{X_{-\tau(k,j)}}(dx)\label{eq132}\ee
is a weakly consistent estimate of the true conditional distribution 
$P(dx|X^-Y^-Y)$. Thus all results in this paper remain valid if the 
decisions can be made with knowledge of not only 
the past but also side information.

\vfill\eject
\centerline{\large\bf Appendix}

\smallskip
Let $\{\B_k\}_{k\ge 1}$ be an increasing 
sequence of finite subfields that asymptotically generate 
the Borel $\sigma$-field on $\X$, and suppose the empirical
conditional distributions $\hat P_k(dx|X^{-\lambda(k)})$ 
are defined as in (\ref{eq16}). Let $T$ denote the left shift 
on the two-sided sequence space $\X^\Z$.

\smallskip
\noindent
{\bf Theorem 1A.} {\sl If $\{X_t\}$ is a stationary process  
with values in a complete separable metric (Polish) space $\X$ 
then for every 
set $B$ in the generating field $\bigcup_k\,\B_k$, we have
\be\lim_k\hat P_k\{X\in B|X^{-\lambda_k}\}
=P\{X\in B|X^-\}\quad\hbox{in $L^1$.}\label{eq133}\ee
}

\noindent
{\sl Proof:} It follows from the martingale convergence
theorem that
\be\lim_k\,P\{X\in B|[X^{-\l_k}]^k\}=P\{X\in B|X^-\}\label{eq134}\ee
almost surely and in $L^1$. Thus it suffices to show that 
$E|\Theta_k|\to 0$ where
\be\Theta_k={1\over J_k}\sum_{0\le j\le J_k}1\{X_{-\tau(k,j)}\in B\}
-P\{X_0\in B|[X^{-\l_k}]^k\}.\label{eq135}\ee
We claim that $E|\Theta_k|=E|\tilde\Theta_k|$ where
\be\tilde \Theta_k={1\over J_k}\sum_{0\le j\le J_k}
1\{X_{\tilde\tau(k,j)}\in B\}-P\{X_0\in B|[X^{-\l_k}]^k\}.\label{eq136}\ee
Indeed, for any measurable function $g(\Theta)\ge 0$ 
(including $g(\Theta)=|\Theta|$) and for any integer
sequence $0=t_0<t_1<\ldots<t_{J(k)}=t$, we have 
\be[1\{\tau^k_j=t_j,\;0\le j\le J_k\}\,g(\Theta_k)]=
[1\{\tilde\tau^k_{J_k-j}=t-t_j,\;0\le j\le J_k\}\,
g(\tilde\Theta_k)]\circ T^{-t}\label{eq137}\ee
and consequently, by stationarity,
\begin{eqnarray}
Eg(\Theta_k)
&=&\sum_{t}\,\sum_{0=t_0<t_1<\ldots<t_{J_k}=t}\,
E\{1\{\tau^k_j=t_j,\;0\le j\le J_k\}\,g(\Theta_k)\}\nonumber \\
&=&\sum_{t}\,\sum_{0=t_0<t_1<\ldots<t_{J_k}=t}\,
E\{1\{\tilde\tau^k_{J_k-j}=t-t_j,\;0\le j\le J_k\}
\,g(\tilde\Theta_k)\} \nonumber \\
&=&\sum_{t}\,\sum_{0=\tilde t_0<\tilde t_1<\ldots<\tilde t_{J_k}=t}\,
E\{1\{\tilde\tau^k_i=\tilde t_i,\;0\le i\le J_k\}
\,g(\tilde\Theta_k)\}=Eg(\tilde \Theta_k). 
\label{eq138}
\end{eqnarray}
Observe that $\{\tilde\tau^k_j-1\}_{j\ge 0}$ is 
an increasing sequence of stopping times adapted 
to the filtration $\{\F^k_t\}_{t\ge 0}$ where
\be\F^k_t=\sigma(\ldots,[X_{-1}]^k,[X_0]^k,[X_1]^k,
\ldots,[X_{t-1}]^k).\label{eq139}\ee
Let $\tilde\F^k_j$ denote the $\sigma$-field of events 
that are expressible in terms of the quantized random 
variables $[X_t]^k$ at times $t<\tilde\tau^k_j$. Thus 
$\tilde\F^k_j$ is the $\sigma$-field of events $F$ such 
that $F\bigcap \{\tilde\tau^k_j=t\}$ belongs to $\F^k_t$ 
for all $t\ge 0$, and $\tilde\F^k_j$ is generated by 
the family of events 
$\{F_t\bigcap\{\tau^k_j=t\}:\;F_t\in\F^k_t,\;t\ge 0\}$.
One may decompose $\tilde\Theta_k$ into the sum
\be\tilde\Theta_k=
{1\over J_k}\sum_{0\le j\le J_k}(\Delta^k_j+\Phi^k_j),\label{eq140}\ee
where
\be\Delta^k_j= 1\{X_{\tilde\tau(k,j)}\in B\}-
P\{X_{\tilde\tau(k,j)}\in B|\tilde\F^k_j\},\label{eq141}\ee
\be\Phi^k_j= P\{X_{\tilde\tau(k,j)}\in B|\tilde\F^k_j\}
-P\{X_0\in B|[X^{-\l_k}]^k\}.\label{eq142}\ee
Notice that $\{\Delta^k_j\}_{j\ge 0}$ is a martingale 
difference sequence with respect to the filtration 
$\{\tilde\F^k_j\}_{j\ge 0}$ (in the sense that 
$\Delta^k_j$ is $\tilde\F^k_{j+1}$-measurable 
and $E\{\Delta^k_j|\tilde\F^k_j\}=0$ for all $j\ge 0$). 
Since $|\Delta^k_j|\le 1$ and the random variables
$\Delta^k_j$ are orthogonal, we see that 
\be E\bigg|{1\over J_k}\sum_{0\le j\le J_k}\Delta^k_j\bigg|^2
={1\over J_k^2}\sum_{0\le j\le J_k}E|\Delta^k_j|^2
\le \left({1+J_k\over J^2_k}\right)\label{eq143}\ee
and consequently (since $(E|Z|)^2\le E\{|Z|^2\}$ and $J_k\to\infty$),
\be E\bigg|{1\over J_k}\sum_{0\le j\le J_k}\Delta_j^k
\bigg|\le {\sqrt{1+J_k}\over J_k}\to 0.\label{eq144}\ee
Also observe that for any measurable function 
$g(\Phi)\ge 0$ and any integer $t\ge 0$, 
\be[1\{\tau^k_j=t\}\,g(\Phi^k_0)]=
[1\{\tilde\tau^k_j=t\}\,g(\Phi^k_j)]\circ T^{-t}\label{eq145}\ee
and consequently, by stationarity,
\begin{eqnarray}
Eg(\Phi^k_j)
&=&\sum_{t\ge 0}E\{1\{\tilde\tau^k_j=t\}\,g(\Phi^k_j)\} \nonumber \\
&=&\sum_{t\ge 0}E\{1\{\tau^k_j=t\}\,g(\Phi^k_0)\}
=Eg(\Phi^k_0).
\label{eq146}
\end{eqnarray}
In particular, setting $g(\Phi)=|\Phi|$ proves 
that $E|\Phi^k_j|=E|\Phi^k_0|$. 
By the martingale convergence theorem 
\be\Phi^k_0=P\{X_0\in B|[X^{-\infty}]^k\}-
P\{X_0\in B|[X^{-\l_k}]^k\}\to 0
\quad\hbox{almost surely and in $L^1$}\label{eq147}\ee
and consequently
\be E\bigg|{1\over J_k}\sum_{0\le j\le J_k}
\Phi^k_j\bigg|\le {1\over J_k}\sum_{0\le j\le J_k}
E|\Phi^k_j|=\left({1+J_k\over J_k}\right)\,E|\Phi^k_0|\to 0.\label{eq148}\ee
The desired conclusion $E|\Theta_k|\to 0$ follows since
\be E|\Theta_k|=E|\tilde\Theta_k|\le
E\bigg|{1\over J_k}\sum_{0\le j\le J_k}\Delta_j^k\bigg|
+E\bigg|{1\over J_k}\sum_{0\le j\le J_k}
\Phi_j^k\bigg|\to 0.\label{eq149}\ee
\qed

\vfill\eject
\bigskip
\centerline{\large\bf Acknowledgment}

\medskip
The authors wish to express their gratitude to L\'{a}szl\'{o} 
Gy\"orfi for his encouragement and suggestions 
regarding this investigation. The first author 
thanks Tam\'{a}s  Szabados for the many discussions 
on Polish spaces. The second author's efforts 
have been partially supported by NSF Grant
INT 92 01430 as well as NIH Grant R01 AI37535.


\begin{thebibliography}{99}

\bibitem{bib1}
P. H. Algoet, 
"Universal schemes for prediction, gambling and 
portfolio selection,"
{\sl Annals Probab.,} vol 20, pp. 901--941, 1992.
Correction: {\sl ibid.}, vol. 23, pp. 474--478, 1995.

\bibitem{bib2}
P. H. Algoet,
"The strong law of large numbers for sequential
decisions under uncertainty," 
{\sl IEEE Trans. Inform. Theory,} 
vol. 40, pp. 609--634, May 1994.

\bibitem{bib3}
P. H. Algoet and T. M. Cover,
"Asympotic optimality and asymptotic equipartition
properties of log-optimum investment,"
{\sl Annals Probab.,} vol. 16, pp. 876--898, 1988.

\bibitem{bib4}
K. Azuma,
"Weighted sums of certain dependent random variables,"
{\sl Tohoku Mathematical Journal,} vol. 37, pp. 357--367, 1967.

\bibitem{bib5}
D. H. Bailey,
{\sl Sequential Schemes for Classifying and Predicting 
Ergodic Processes.} Ph. D. thesis, Stanford University, 1976.

\bibitem{bib6}
A. R. Barron,
"Entropy and the central limit theorem,"
{\sl Ann. Probab.,} vol. 14, pp. 336--342, 1986.

\bibitem{bib7}
A. R. Barron, L. Gy\"orfi and E. C. van der Meulen,
"Distribution estimation consistent in total variation and
in two types of information divergence,"
{\sl IEEE Trans. Inform. Theory,} 
vol. 38, pp. 1437--1454, Sept. 1992.

\bibitem{bib8}
L. Breiman,
"The individual ergodic theorem of information theory,"
{\sl Ann. Math. Statist.,} vol. 28, pp. 809--811, 1957.
Correction: {\sl ibid.,} vol. 31, pp. 809--810, 1960.

\bibitem{bib9}
T. M. Cover,
"Open problems in information theory,"
in {\sl 1975 IEEE Joint Workshop on Information Theory,}
pp. 35--36. New York: IEEE Press, 1975.

\bibitem{bib10}
T. M. Cover and J. Thomas,
{\sl Elements of Information Theory.}
New York: Wiley, 1991.

\bibitem{bib11}
I. Csisz\'ar and J. K\"orner, 
{\sl Information Theory: Coding Theorems 
for Discrete Memoryless Systems.}
Budapest: Akad\'emiai Kiad\'o, 1981.

\bibitem{bibFeder} 
M. Feder, N. Merhav and M. Gutman,
"Universal prediction of individual sequences,"
{\sl IEEE Trans. Inform. Theory,} vol. 38, pp. 1258--1270, July 1992. 

\bibitem{bib13}
A. Gavish and A. Lempel,
"Match-length functionals for data compression,"
presented at {\sl IEEE Int. Symp. Inform. Theory,} 
Trondheim, Norway, Jun. 27--Jul. 1, 1994.

\bibitem{bib14} 
R. M. Gray, 
{\sl Probability, Random Processes, and Ergodic Theory.}
New York: Springer-Verlag, 1988.

\bibitem{bib15}
L. Gy\"orfi, W. H\"ardle, P. Sarda and Ph. Vieu,
{\sl Nonparametric Curve Estimation from Time Series.}
Berlin: Springer-Verlag, 1989.

\bibitem{bib16}
L. Gy\"orfi, I. P\'ali and E. C. van der Meulen,
"There is no universal source code for infinite source alphabet,"
{\sl IEEE Trans. Inform. Theory,} vol. 40, pp. 267--271, Jan. 1994.

\bibitem{bib17}
M. Kac,
"On the notion of recurrence in discrete stochastic processes,"
{\sl Bull. Amer. Math. Soc.}, vol. 53, pp. 1002--1010, Oct. 1947.



\bibitem{bib18}
G. G. Langdon, Jr.,
"A note on the Lempel-Ziv model for compressing
individual sequences,"
{\sl IEEE Trans. Inform. Theory,} 
vol. IT-29, pp. 284--287, Mar. 1994.

\bibitem{bib19}
K. Marton and P. C. Shields,
"Entropy and the consistent estimation of
joint distributions,"
{\sl Annals Probab.,} vol. 22, pp. 960--977, Apr. 1994. 

\bibitem{bib20}
G. Morvai, 
{\sl Estimation of Conditional Distributions for
Stationary Time Series.}
Ph. D. Thesis, Technical University of Budapest, 1994.

\bibitem{bib21}
G. Morvai, S. Yakowitz, and L. Gy\"orfi,
"Nonparametric inferences for ergodic, stationary time series,"
 {\sl The Annals of Statistics,} vol. 24, pp. 370-379, 1996.

\bibitem{bib22}
D. S. Ornstein,
"Guessing the next output of a stationary process,"
{\sl Israel J. Math.,} vol. 30, pp. 292--296, 1978.

\bibitem{bib23}
D. S. Ornstein and P. C. Shields,
"Universal almost sure data compression,"
{\sl Annals Probab.}, vol. 18, pp. 441--452, 1990.


\bibitem{bib24}
D. S. Ornstein and B. Weiss,
"Entropy and data compression schemes,"
{\sl IEEE Trans. Inform. Theory,} vol. 39, pp. 78--83, Jan. 1993.

\bibitem{bib25}
A. Perez, 
"On Shannon-McMillan's limit theorem for pairs 
of stationary processes,"
{\sl Kybernetika,} vol. 16, nr. 4, pp. 301--314, 1980.

\bibitem{bib26}
M. S. Pinsker,
{\sl Information and Information Stability 
of Random Variables and Processes.} 
Translated and edited by A. Feinstein.
San Francisco: Holden-Day, 1964.

\bibitem{bib27}
J. Rissanen,
"Stochastic complexity and modeling,"
{\sl Ann. Statist.,} vol. 14, pp. 1080--1100, 1986.

\bibitem{bib28}
J. Rissanen and G. G. Langdon, Jr.,
"Universal modeling and coding,"
{\sl IEEE Trans. Inform. Theory}, 
vol. IT-27, pp. 12--23, Jan. 1981.



\bibitem{bib29}
B. Ya. Ryabco,
"Twice-universal coding,"
{\sl Problems of Inform. Trans.,} 
vol. 20, pp. 173--177, July-Sept. 1984.


\bibitem{bib30}
B. Ya. Ryabco,
"Prediction of random sequences and universal coding,"
{\sl Problems of Inform. Trans.,} 
vol. 24, pp. 87-96, Apr.-June 1988.


\bibitem{bib32}
B. Scarpellini,
"Conditional expectations of stationary processes,"
{\sl Z. Wahrscheinlichkeitstheorie verw. Gebiete,}
vol. 56, pp. 427--441, 1981.

\bibitem{bib34}
P. C. Shields, 
"Universal redundancy rates don't exist,"
{\sl IEEE Trans. Inform. Theory}, vol. 39, pp. 520--524, Mar. 1993.

\bibitem{bib35}
F. M. J. Willems, 
"Universal data compression and repetition times,"
{\sl IEEE Trans. Inform. Theory,} vol. 35, pp. 54--58, Jan. 1989.

\bibitem{bib36}
A. D. Wyner and J. Ziv,
"Some asymptotic properties of entropy of a stationary
ergodic data source with applications to data compression,"
{\sl IEEE Trans. Inform. Theory,} vol. 35, pp. 1250--1258, Nov. 1989.

\bibitem{bib37}
S. Yakowitz, L. Gy\"orfi, and G. Morvai,
"An algorithm for nonparametric forecasting 
for ergodic, stationary time series,"
presented at {\sl IEEE Int. Symp. Inform. Theory,} 
Trondheim, Norway, Jun. 27--Jul. 1, 1994.

\bibitem{bib38}
J. Ziv and A. Lempel,
"Compression of individual sequences by 
variable rate coding,"
{\sl IEEE Trans. Inform. Theory,} vol. IT-24, 
pp. 530--536, Sept. 1978.
 
\end{thebibliography}
\end{document}